%% file: admm.tex
\documentclass[10pt,twocolumn,twoside]{IEEEtran}
\usepackage[bookmarks=false,hyperfigures=false]{hyperref}
\usepackage{amsmath,amssymb,mathtools,bbm,url,doi,xspace,relsize}
\usepackage{subfig}

\pdfoutput=1
\usepackage{ifpdf}
\ifpdf
  \usepackage{graphicx,color}
  \DeclareGraphicsExtensions{.pdf}
\else
  \usepackage{graphicx,color}
  \DeclareGraphicsExtensions{.pstex,.ps,.eps}
\fi

\newcommand{\ie}{{i.e.}\xspace}
\newcommand{\eg}{{e.g.}\xspace}
\newcommand{\etal}{{et al.}\xspace}
\newcommand{\eq}[1]{Eq.~(\ref{eq:#1})}
\newcommand{\eqc}[1]{(\ref{eq:#1})}
\newcommand{\fig}[1]{Fig.~\ref{fig:#1}}
\newcommand{\figs}[1]{Figs.~\ref{fig:#1}}
\newcommand{\fign}[1]{\ref{fig:#1}}

\newcommand{\sctn}[1]{Sec.~\ref{sec:#1}}

\newcommand{\tl}[1]{\tilde{#1}}
\newcommand{\mb}[1]{\mathbf{#1}}
\newcommand{\tmb}[1]{\tilde{\mathbf{#1}}}

\newcommand{\mbb}[1]{\mathbb{#1}}

\newcommand{\norm}[1]{\left\| #1 \right\|}
\DeclarePairedDelimiterX{\normsz}[1]{\lVert}{\rVert}{#1}

\DeclareMathOperator*{\argmin}{arg\,min}

\DeclareRobustCommand{\bvsigma}{\ensuremath{\boldsymbol{\varsigma}}}

\newcommand{\resizemath}[2]{%
  \resizebox{#1}{!}{$ \displaystyle
  #2
  $}
}%

\makeatletter%
\if@twocolumn%
\else%
\renewcommand{\resizemath}[2]{%
  #2
}%
\fi
\makeatother

\newcommand{\smallmath}[1]{%
  {\small
    \setlength{\abovedisplayskip}{6pt}
    \setlength{\belowdisplayskip}{\abovedisplayskip}
    \setlength{\abovedisplayshortskip}{0pt}
    \setlength{\belowdisplayshortskip}{3pt}
    #1
  }%
}%

\begin{document}

\title{ADMM Penalty Parameter Selection by Residual Balancing}

\author{Brendt Wohlberg
  \thanks{Theoretical Division, Los Alamos
    National Laboratory, Los Alamos, NM 87545, USA. Email:
    \texttt{brendt@lanl.gov}, Tel: +1 505 667 6886
  }
  \thanks{This research was supported by the U.S. Department
    of Energy through the LANL/LDRD Program.}
}

\maketitle

\begin{abstract}
  Appropriate selection of the penalty parameter is crucial to
  obtaining good performance from the Alternating Direction Method of
  Multipliers (ADMM). While analytic results for optimal selection of
  this parameter are very limited, there is a heuristic method that
  appears to be relatively successful in a number of different
  problems. The contribution of this paper is to demonstrate that
  their is a potentially serious flaw in this heuristic approach, and
  to propose a modification that at least partially addresses it.
\end{abstract}


\IEEEpeerreviewmaketitle

\section{Introduction}
\label{sec:intro}

The Alternating Direction Method of Multipliers (ADMM) has become a
very popular approach to solving a broad variety of optimization
problems in signal and image processing, prominent examples including
Total Variation regularization and sparse representation
problems~\cite{goldstein-2009-split},~\cite[Sec.
6]{boyd-2010-distributed}, ~\cite{afonso-2011-augmented}. This method
introduces an additional parameter, the \emph{penalty parameter}, on
which the rate of convergence is strongly dependent, but for which
there are no analytic results to guide selection other than for a very
specific set of problems~\cite{ghadimi-2015-optimal,
  raghunathan-2014-alternating},~\cite[Sec. 5]{raghunathan-2015-admm}. There
is, however, a heuristic method for automatically adapting the penalty
parameter~\cite{he-2000-alternating} that appears to becoming quite
popular~\cite{hansson-2012-subspace, liu-2013-nuclear,
  vu-2013-fantope, iordache-2014-collaborative, weller-2014-phase,
  wohlberg-2014-efficient}. The present paper demonstrates a serious
flaw in this heuristic approach, and proposes a modification that at
least partially addresses it.

\section{ADMM}
\label{sec:admm_dtl}

The notation and exposition in this section follows that of the
influential tutorial by Boyd \etal~\cite{boyd-2010-distributed}. The
\emph{Lagrangian} for the constrained problem
\begin{equation}
\argmin_{\mb{x}} f(\mb{x}) \; \text{ such that } \; A \mb{x} = \mb{b}
\;,
\label{eq:linconprimal}
\end{equation}
is
\begin{equation}
L(\mb{x}, \mb{y}) = f(\mb{x}) + \mb{y}^T (A \mb{x} - \mb{b})
\;,
\label{eq:pcnstlgrng}
\end{equation}
where $\mb{x}$ and $\mb{y}$ are referred to as the \emph{primal} and
\emph{dual} variables respectively. The primal and dual feasibility
conditions
\begin{align}
0 = \nabla L(\mb{x}^*, \cdot)
\;\;\Rightarrow\;\; & A \mb{x}^* - \mb{b} = 0  \label{eq:lgrngnopty}\\
0 \in \partial L(\cdot, \mb{y}^*)
\;\;\Rightarrow\;\; & 0 \in \partial f(\mb{x}^*)
+ A^T \mb{y}^* \;, \label{eq:lgrngoptx}
\end{align}
where $\partial$ denotes the subdifferential
operator~\cite[Ch. D]{urruty-2004-fundamentals}, provide conditions on
the optimal primal and dual variables $\mb{x}^*$ and $\mb{y}^*$.
The \emph{method of multipliers} solves this problem
via \emph{dual ascent}
\begin{align}
\mb{x}^{(k+1)} & = \argmin_{\mb{x}} L_{\rho}(\mb{x}, \mb{y}^{(k)}) \\
\mb{y}^{(k+1)} & = \mb{y}^{(k)} + \rho (A \mb{x}^{(k+1)} -
\mb{b}) \; ,
\end{align}
where $L_{\rho}$ is the \emph{augmented Lagrangian}
\begin{equation}
L_{\rho}(\mb{x}, \mb{y}) = f(\mb{x}) + \mb{y}^T (A \mb{x} - \mb{b}) +
\frac{\rho}{2} \norm{A \mb{x} - \mb{b}}_2^2  %
\end{equation}
with \emph{penalty parameter} $\rho$.

ADMM can be viewed as a variant of this method\footnote{There are limitations to this interpretation~\cite{eckstein-2012-augmented}.} applied to the problem
\begin{equation}
\argmin_{\mb{x},\mb{z}} f(\mb{x}) + g(\mb{z}) \; \text{ such that } \;
A \mb{x} + B \mb{z} = \mb{c} \; ,
\label{eq:admmprob}
\end{equation}
where $\mb{x} \in \mbb{R}^n$, $\mb{z} \in \mbb{R}^m$, and $\mb{c} \in
\mbb{R}^p$, and the Lagrangian and augmented Lagrangian are, respectively,
\begin{align}
L(\mb{x}, \mb{z}, \mb{y}) = & f(\mb{x}) + g(\mb{z}) + \mb{y}^T
(A \mb{x} + B \mb{z} - \mb{c}) \\
L_{\rho}(\mb{x}, \mb{z}, \mb{y}) = & L(\mb{x}, \mb{z}, \mb{y}) +
\frac{\rho}{2} \norm{A \mb{x} + B \mb{z} - \mb{c}}_2^2  \; .
\end{align}
Instead of jointly solving for $\mb{x}$ and
$\mb{z}$, ADMM  alternates the $\mb{x}$ and
$\mb{z}$ updates  (thus the
\emph{alternating direction})
\begin{align}
\mb{x}^{(k+1)} & = \argmin_{\mb{x}} L_{\rho}(\mb{x}, \mb{z}^{(k)},
\mb{y}^{(k)}) \label{eq:admmx} \\
\mb{z}^{(k+1)} & = \argmin_{\mb{z}} L_{\rho}(\mb{x}^{(k+1)}, \mb{z},
\mb{y}^{(k)}) \label{eq:admmz}  \\
\mb{y}^{(k+1)} & = \mb{y}^{(k)} + \rho (A \mb{x}^{(k+1)} + B
\mb{z}^{(k+1)} - \mb{c}) \; . \label{eq:admmy}
\end{align}

It is often more convenient to work with the \emph{scaled form} of
ADMM, which is obtained by the change of variable to the \emph{scaled
  dual variable} $\mb{u} = \rho^{-1} \mb{y}$. Defining the residual
\begin{equation}
\mb{r} = A \mb{x} + B \mb{z} - \mb{c}
\end{equation}
and replacing $\mb{y}$ with
$\mb{u}$ we have
 \begin{align}
 L_{\rho}(\mb{x}, \mb{z}, \mb{u}) &=  f(\mb{x}) + g(\mb{z}) +
 \frac{\rho}{2} \norm{\mb{r}+\mb{u}}_2^2 - \frac{\rho}{2}
 \norm{\mb{u}}_2^2
 \; .
 \end{align}
Since the minimisers of $L_{\rho}(\mb{x}, \mb{z}, \mb{u})$ with respect
to $\mb{x}$ and $\mb{z}$ do not depend on the final $\frac{\rho}{2}
\norm{\mb{u}}_2^2$ term, %
the iterations can be written as
\smallmath{
\begin{align}
\mb{x}^{(k+1)} & = \argmin_{\mb{x}} f(\mb{x}) + \frac{\rho}{2} \norm{A
  \mb{x} + B \mb{z}^{(k)} - \mb{c} + \mb{u}^{(k)}}_2^2
\label{eq:admmscaledx} \\
\mb{z}^{(k+1)} & = \argmin_{\mb{z}} g(\mb{z}) + \frac{\rho}{2} \norm{A
  \mb{x}^{(k+1)} + B \mb{z} - \mb{c} + \mb{u}^{(k)}}_2^2  \\
\mb{u}^{(k+1)} & =  \mb{u}^{(k)} + A \mb{x}^{(k+1)} + B  \mb{z}^{(k+1)}
- \mb{c}   \; .
\end{align}
}

\subsection{ADMM Residuals}
\label{sec:resid}

Denote optimal primal variables by $\mb{x}^*$ and $\mb{z}^*$, and the
optimal dual variable by $\mb{y}^*$. It will also be useful to define
$p^* = f(\mb{x}^*) + g(\mb{z}^*)$ and $p^{(k)} = f(\mb{x}^{(k)}) +
g(\mb{z}^{(k)})$.  The primal feasibility condition
\begin{equation}
A \mb{x}^* + B \mb{z}^* - \mb{c} =
0 \;,  \label{eq:admmlgprim}
\end{equation}
and dual feasibility conditions
\begin{align}
0 \in \partial L(\cdot, \mb{z}^*, \mb{y}^*)
  \; \Rightarrow \;\; & 0 \in \partial f(\mb{x}^*) + A^T
  \mb{y}^*   \label{eq:admmlgrduafsg} \\
0 \in \partial L(\mb{x}^*, \cdot, \mb{y}^*)
 \; \Rightarrow \;\; & 0 \in \partial g(\mb{z}^*) + B^T
 \mb{y}^*   \label{eq:admmlgrduagsg}
\end{align}
for~\eq{admmprob} hold at the problem solution
$(\mb{x}^*,\mb{z}^*,\mb{y}^*)$. These conditions can be used to derive
convergence measures for ADMM algorithm iterates
$(\mb{x}^{(k)},\mb{z}^{(k)},\mb{y}^{(k)})$.

A natural measure of primal feasibility based on~\eq{admmlgprim} is
the \emph{primal residual}
\begin{equation}
  \mb{r}^{(k+1)} = A \mb{x}^{(k+1)} + B \mb{z}^{(k+1)} - \mb{c} \;.
\label{eq:prires}
\end{equation}
Now, since $\mb{z}^{(k+1)}$ minimises $L_{\rho}(\mb{x}^{(k+1)}, \mb{z},
\mb{y}^{(k)})$ (see~\eq{admmz}), we have
\begin{align}
0 & \in [\partial L_{\rho}(\mb{x}^{(k+1)}, \cdot, \mb{y}^{(k)})](\mb{z}^{(k+1)}) \nonumber \\
   &  = \partial g(\mb{z}^{(k+1)}) + B^T \mb{y}^{(k)}  \nonumber \\
   &  \qquad \qquad \;\;\;\;\;\;\, + \rho B^T (A
   \mb{x}^{(k+1)} + B \mb{z}^{(k+1)} - \mb{c}) \nonumber \\
   & = \partial g(\mb{z}^{(k+1)}) + B^T \mb{y}^{(k)} + \rho B^T
   \mb{r}^{(k+1)} \nonumber \\
   & = \partial g(\mb{z}^{(k+1)}) + B^T (\mb{y}^{(k)} + \rho
   \mb{r}^{(k+1)}) \nonumber \\
   & = \partial g(\mb{z}^{(k+1)}) + B^T \mb{y}^{(k+1)} \;,
 \end{align}
 so that iterates $\mb{z}^{(k+1)}$ and $\mb{y}^{(k+1)}$ always satisfy
 dual feasibility condition~\eq{admmlgrduagsg},
 leaving~\eq{admmlgrduafsg} as the remaining optimality criteria to be
 satisfied.
 Following a similar derivation, since $\mb{x}^{(k+1)}$ minimises
 $L_{\rho}(\mb{x}, \mb{z}^{(k)}, \mb{y}^{(k)})$ (see~\eq{admmx}), we
 have
\begin{align}
0 & \in [\partial L_{\rho}(\cdot, \mb{z}^{(k)}, \mb{y}^{(k)})](\mb{x}^{(k+1)}) \nonumber \\
  & = \partial f(\mb{x}^{(k+1)}) + A^T \mb{y}^{(k)}  + \rho A^T (A
   \mb{x}^{(k+1)} + B \mb{z}^{(k)} - \mb{c}) \nonumber \\
  & = \partial f(\mb{x}^{(k+1)}) + A^T \mb{y}^{(k)} + \rho A^T (A
  \mb{x}^{(k+1)} + B
 \mb{z}^{(k+1)} \nonumber \\
 &  \qquad \qquad \qquad \qquad \qquad \;\;\;\;\; - \mb{c} + B
 \mb{z}^{(k)} - B  \mb{z}^{(k+1)}) \nonumber\\
  & = \partial f(\mb{x}^{(k+1)}) + A^T \mb{y}^{(k)} \nonumber \\
&  \qquad \qquad \qquad  + \rho A^T
  (\mb{r}^{(k+1)} + B \mb{z}^{(k)} - B  \mb{z}^{(k+1)}) \nonumber \\
 & = \partial f(\mb{x}^{(k+1)}) + A^T (\mb{y}^{(k)} \nonumber \\
&  \qquad \qquad \qquad  + \rho
 \mb{r}^{(k+1)}) + \rho A^T B ( \mb{z}^{(k)} -   \mb{z}^{(k+1)})\nonumber \\
 & = \partial f(\mb{x}^{(k+1)}) + A^T \mb{y}^{(k+1)} \!+\! \rho A^T B (
 \mb{z}^{(k)} \!-\! \mb{z}^{(k+1)})  \;.
\label{eq:dualres0}
\end{align}
Setting $\rho A^T B (\mb{z}^{(k+1)} -
\mb{z}^{(k)}) = 0$ in~\eq{dualres0} implies that $\mb{x}^{(k+1)},
\mb{y}^{(k+1)}$ satisfy dual feasibility condition~\eq{admmlgrduafsg},
which suggests defining
\begin{equation}
  \mb{s}^{(k+1)} = \rho A^T B (\mb{z}^{(k+1)} - \mb{z}^{(k)})
\label{eq:dualres}
\end{equation}
as the \emph{dual residual} based on dual feasibility
condition~\eq{admmlgrduafsg}.

Since both primal and dual residuals converge to zero as the ADMM
algorithm progresses~\cite[Sec. 3.3]{boyd-2010-distributed}, they can
be used to define ADMM algorithm convergence measures. It is also
worth noting that~\eq{admmx} and~\eq{admmz} suggest that the norm of
the primal residual decreases with increasing $\rho$ (and vice versa),
and the definition of the dual residual suggests that it increases
with increasing $\rho$ (and vice versa).

\subsection{Adaptive Penalty Parameter}
\label{sec:adaptrho}

As discussed in~\sctn{intro}, the correct choice of the penalty
parameter plays a vital role in obtaining good convergence.  He
\etal~\cite{he-2000-alternating} define the distance from convergence
as $\| \mb{r}^{(k+1)} \|_2^2 + \| \mb{s}^{(k+1)}\|_2^2$, and argue
that adaptively choosing the penalty parameter to balance these two
terms is a reasonable heuristic for minimising this distance. This
heuristic is implemented as the update scheme
\begin{equation}
\rho^{(k+1)} = \left\{ \
\begin{array}{ll}
  \tau \rho^{(k)}      & \text{ if } \normsz[\big]{\mb{r}^{(k)}}_2 > \mu
                         \normsz[\big]{\mb{s}^{(k)}}_2\\[3pt]
\tau^{-1} \rho^{(k)}  & \text{ if } \normsz[\big]{\mb{s}^{(k)}}_2 > \mu
                         \normsz[\big]{\mb{r}^{(k)}}_2\\[3pt]
\rho^{(k)}           & \text{ otherwise } \;,
\end{array}
\right.
\label{eq:rhoupdate}
\end{equation}
where $\tau$ and $\mu$ are constants, the usual values being $\tau =
2$ and $\mu = 10$~\cite{he-2000-alternating,
  wang-2001-decomposition},~\cite[Sec 3.4.1]{boyd-2010-distributed}.

This scheme has has been found to be effective
for a variety of problems~\cite{hansson-2012-subspace,
  liu-2013-nuclear, vu-2013-fantope, iordache-2014-collaborative,
  weller-2014-phase, wohlberg-2014-efficient}, but
it will be demonstrated that it suffers from a potentially serious flaw.

\subsection{Stopping Criteria}
\label{sec:stopcrit}

The residuals can be used to define stopping criteria for the ADMM
iterations; \eg Boyd \etal~\cite[Sec. 3.3.1]{boyd-2010-distributed}
recommend stopping criteria
\begin{equation}
  \normsz[\big]{\mb{r}^{(k)}}_2 \leq \epsilon_{\mathrm{pri}}^{(k)}
  \;\; \text{ and }
  \normsz[\big]{\mb{s}^{(k)}}_2 \leq \epsilon_{\mathrm{dua}}^{(k)}
\end{equation}
where
\smallmath{
\begin{align}
 \epsilon_{\mathrm{pri}}^{(k)} & \!=\! \sqrt{p} \epsilon_{\mathrm{abs}} \!+\!
 \epsilon_{\mathrm{rel}} \max\left\{\normsz[\big]{A \mb{x}^{(k)}}_2,
   \normsz[\big]{B
      \mb{z}^{(k)}}_2, \normsz[\big]{\mb{c}}_2 \right\} \label{eq:epri} \\
   \epsilon_{\mathrm{dua}}^{(k)} \!&\! = \sqrt{n} \epsilon_{\mathrm{abs}} +
   \epsilon_{\mathrm{rel}} \normsz[\big]{A^T \mb{y}^{(k)}}_2 \label{eq:edua}
   \;,
\end{align}}
$\epsilon_{\mathrm{abs}}$ and $\epsilon_{\mathrm{rel}}$ are
absolute and relative tolerances respectively, and $n$ and $p$ are the
dimensionalities of $\mb{x}$ and $\mb{c}$ respectively (\ie $\mb{x}
\in \mbb{R}^n$ and $\mb{c} \in \mbb{R}^p$).

\section{ADMM Problem Scaling Properties}
\label{sec:admmscale}

Let us consider the behaviour of ADMM under scaling of the
optimization problem being addressed, denoting~\eq{admmprob} as
problem $P$, and defining $\tilde{P}$ as
\begin{equation}
  \argmin_{\mb{x},{\mb{z}}} \alpha f(\gamma
  \mb{x}) + \alpha g(\gamma \mb{z}) \; \text{ s.t. } \;
  \beta A \gamma \mb{x} + \beta B \gamma \mb{z} =
  \beta \mb{c}  \;.
\label{eq:admmp1}
\end{equation}
In this problem $\alpha$ represents a scaling of the objective
function, $\beta$ represents a scaling of the constraint, and $\gamma$
represents a scaling of the problem variables. These scalings are
chosen to parameterise the family of scalings of an ADMM problem under
which the solution is invariant, modulo a scaling\footnote{The
  minimisers of $\tilde{P}$ are invariant to $\alpha$ and $\beta$, and
  are invariant to $\gamma$ modulo a scaling factor.}. It is important to
emphasise that these scalings can represent both explicit scaling of a
problem and the implicit scaling with respect to alternative possible
choices\footnote{For problems involving physical quantities, for
  example, scaling by $\alpha$ and $\gamma$ correspond respectively to
  choices of the units in which the functional value and solution are
  expressed.  Scaling by $\beta$ corresponds to the choices to be made
  in constructing the constraint; for example, if $\mb{z}$ is to
  represent the gradient of $\mb{x}$, then $A$ could be scaled to
  represent differences between samples with or without normalisation
  by the physical step size of the grid on which $\mb{x}$ is defined.}
inherent in choosing functional, constraints, and variables.
Problem $\tilde{P}$ can be expressed in the standard form as
\begin{equation}
  \argmin_{\mb{x},{\mb{z}}} \tl{f}(
  \mb{x}) + \tl{g}(\mb{z}) \; \text{ such that } \;
  \tl{A} \mb{x} + \tl{B} \mb{z} = \tmb{c}
\label{eq:admmp1std}
\end{equation}
with
\begin{gather}
  \tl{f}(\mb{x}) = \alpha f(\gamma \mb{x}) \quad
  \tl{g}(\mb{z}) = \alpha g(\gamma \mb{z}) \quad \nonumber \\
  \tl{A} = \beta \gamma A \quad
  \tl{B} = \beta \gamma B \quad
  \tmb{c} = \beta \mb{c} \;.
\end{gather}

The Lagrangian is
\begin{align}
  \tilde{L}({\mb{x}}, {\mb{z}}, {\mb{y}}) =
  \alpha f(\gamma \mb{x}) & + \alpha g(\gamma
  \mb{z}) \nonumber \\ & + \mb{y}^T (\beta \gamma A \mb{x} +
   \beta \gamma B \mb{z} -  \beta \mb{c})
\; ,
\end{align}
and the primal and dual feasibility conditions are
\begin{equation}
  \beta A \gamma \tmb{x}^* + \beta B \gamma \tmb{z}^* -
  \beta \tmb{c} = 0 \;,
\end{equation}
and
\smallmath{
\begin{align}
  0 \in \partial \tilde{L}(\cdot, \tilde{\mb{z}}^*, \tilde{\mb{y}}^*)
  \; \Rightarrow \; & 0\!\in\! \alpha \gamma [\partial f(\cdot)](\gamma
  \tilde{\mb{x}}^*)
  + \beta \gamma A^T \tilde{\mb{y}}^* = 0   \\
  0 \in \partial \tilde{L}(\tilde{\mb{x}}^*, \cdot, \tilde{\mb{y}}^*)
 \; \Rightarrow \; & 0\! \in\! \alpha \gamma [\partial g(\cdot)](\gamma
  \tilde{\mb{z}}^*) + \beta \gamma B^T \tilde{\mb{y}}^* = 0
\end{align}
}
respectively.
It is easily verified that if $\mb{x}^*$, $\mb{z}^*$, and $\mb{y}^*$
satisfy the optimality criteria \eq{admmlgprim}, \eqc{admmlgrduafsg},
and \eqc{admmlgrduagsg} for problem $P$, then
\begin{gather}
\tilde{\mb{x}}^* = \gamma ^{-1} \mb{x}^* \quad \;\;
\tilde{\mb{z}}^* = \gamma ^{-1} \mb{z}^* \quad \;\;
\tilde{\mb{y}}^* = \frac{\alpha}{\beta} \mb{y}^*
\label{eq:optscaling}
\end{gather}
satisfy the primal and dual feasibility criteria for $\tilde{P}$.
The augmented Lagrangian for $\tilde{P}$ is
\smallmath{
\begin{align}
  \tilde{L}_{\tilde{\rho}}(\mb{x}, \mb{z}, \mb{y})
  & =
\alpha f(\gamma
  \mb{x}) + \alpha g(\gamma \mb{z}) \nonumber \\ &+ \alpha
  \left(\frac{\beta}{\alpha} \mb{y}^T \right) ( \gamma A \mb{x} +
  \gamma B \mb{z} - \mb{c} ) \nonumber \\ &+ \alpha \left(\frac{\beta^2}{\alpha}
    \tilde{\rho} \right) \frac{1}{2} \norm{ \gamma A \mb{x} +
    \gamma B \mb{z} - \mb{c}}_2^2 \;,
\end{align}
}
so that setting %
\begin{equation}
 \tilde{\rho} = \frac{\alpha}{\beta^2} \rho
\label{eq:rhoscale}
\end{equation}
gives
\begin{equation}
  \tilde{L}_{\tilde{\rho}}(\mb{x}, \mb{z}, \mb{y}) = \alpha
  L_{\rho}\left(\gamma \mb{x}, \gamma \mb{z}, \frac{\beta}{\alpha} \mb{y}\right)
  \;.
\label{eq:tlscale}
\end{equation}

The iterates $\mb{x}^{(k+1)}$, $\mb{z}^{(k+1)}$, and $\mb{y}^{(k+1)}$
for iteration $k$ of the ADMM algorithm for $P$ are given by
\eq{admmx}, \eqc{admmz}, and \eqc{admmy}. We now consider the
corresponding iterates for $\tilde{P}$, assuming that
\begin{gather}
 \tilde{\mb{z}}^{(k)} = \gamma^{-1} \mb{z}^{(k)} \quad \quad
 \tilde{\mb{y}}^{(k)} =  \frac{\alpha}{\beta} \mb{y}^{(k)} \;.
\label{eq:initscaled}
\end{gather}
The $\mb{x}$ update is
\begin{align}
 \tmb{x}^{(k+1)} & = \argmin_{\mb{x}} \tl{L}_{\tl{\rho}}(\mb{x}, \tmb{z}^{(k)},
\tmb{y}^{(k)})  \nonumber \\
  &= \argmin_{\mb{x}} \tl{L}_{\tl{\rho}}(\mb{x}, \gamma^{-1} \mb{z}^{(k)},
\frac{\alpha}{\beta} \mb{y}^{(k)}) \nonumber \\
 &= \argmin_{\mb{x}}  \alpha L_{\rho} (\gamma \mb{x}, \mb{z}^{(k)},
 \mb{y}^{(k)}) \;.
\end{align}
For convex $f$ we have that if $\mb{x}^*$ minimises $f(\mb{x})$ then
$\gamma^{-1} \mb{x}^*$ minimises $\tl{f}(\mb{x}) = \alpha f (\gamma
\mb{x})$, so
\begin{equation}
  \tmb{x}^{(k+1)} = \gamma^{-1} \mb{x}^{(k+1)} \;,
\end{equation}
and similarly it can be shown that
\begin{equation}
  \tmb{z}^{(k+1)} = \gamma^{-1} \mb{z}^{(k+1)} \;.
\end{equation}
For the $\mb{y}$ update we have
\begin{align}
  \tmb{y}^{(k+1)} & = \tmb{y}^{(k)} + \tl{\rho} (\beta \gamma A
  \tmb{x}^{(k+1)} + \beta \gamma B \tmb{z}^{(k+1)} - \beta \mb{c})
  \nonumber \\
  &= \frac{\alpha}{\beta} \left( \mb{y}^{(k)} + \rho (A \mb{x}^{(k+1)}
    + B \mb{z}^{(k+1)} - \mb{c}) \right) \nonumber \\
  &= \frac{\alpha}{\beta} \mb{y}^{(k+1)} \;.
\end{align}

Finally, the primal and dual residuals for $\tilde{P}$ have the
following scaling relationship with those of $P$:
\begin{align}
  \tilde{\mb{r}}^{(k+1)} & = \tl{A} \tmb{x}^{(k+1)} + \tl{B}
  \tmb{z}^{(k+1)} - \tmb{c} \nonumber \\
  &= \beta A \mb{x}^{(k+1)} +
  \beta B \mb{z} ^{(k+1)} - \beta \mb{c} \nonumber \\
  & = \beta\mb{r}^{(k+1)} \label{eq:priresscl} \\
  \tilde{\mb{s}}^{(k+1)} & = \tilde{\rho} \tilde{A}^T \tilde{B}
  (\tilde{\mb{z}}^{(k+1)} - \tilde{\mb{z}}^{(k)}) \nonumber \\
  & = \alpha \gamma \rho A^T B  (\mb{z}^{(k+1)}
  -  \mb{z}^{(k)}) \nonumber \\
  & = \alpha \gamma \mb{s}^{(k+1)} \;. \label{eq:duaresscl}
\end{align}

In summary, the parameters $\alpha$, $\beta$, and $\gamma$ in problem
$\tilde{P}$ generate families of ADMM problems with the same solutions
(modulo a scaling, in the case of $\gamma$), as expressed
in~\eq{optscaling}, but the iterates of the corresponding ADMM
algorithms are only similarly invariant if the initial iterates
(see~\eq{initscaled}) and constant penalty parameter
(see~\eq{rhoscale}) are appropriately scaled.

\section{Residual Balancing}

The scaling properties described in the previous section have a major impact on the residuals and their use within the residual balancing scheme for penalty parameter selection.

\subsection{Adaptive Penalty Parameter}
\label{sec:adaptpenparam}

It was demonstrated above that ADMM algorithm iterates can be made
invariant to problem scaling by a suitable choice of fixed penalty
parameter. It is easily verified that invariance can be maintained
with a varying penalty parameter $\rho^{(k)}$ as long as the required
relationship is also maintained, i.e.
$\tilde{\rho}^{(k)} = \frac{\alpha}{\beta^2} \rho^{(k)}$. If an
adaptive update rule such as~\eq{rhoupdate}, that operates by
multiplying the penalty parameter by some factor, is to preserve this
relationship, it is necessary that (i)
$\tilde{\rho}^{(0)} = \frac{\alpha}{\beta^2} \rho^{(0)}$, and (ii) the
choice of multiplier and when to apply it must be invariant to problem
scaling. But it is clear from~\eq{priresscl} and~\eq{priresscl} that
the primal and dual residuals do not share the same scaling factors,
so that the update rule~\eq{rhoupdate} based on these residuals does
not, in general, preserve the scaling behaviour of the penalty
parameter required to maintain invariance of the algorithm iterates.
It follows that if the adaptive penalty parameter method
of~\sctn{adaptrho} performs well for some problem $P$, it should not
be expected to do so for problem $\tilde{P}$ as the scaling parameters
$\alpha$, $\beta$, and $\gamma$ deviate from unity.

If update rule~\eq{rhoupdate} is known to provide good performance for
 a reference problem $P$, and it becomes necessary to modify
the problem formulation in a way that corresponds to switching to a
scaled problem $\tl{P}$ (e.g. a change of physical units), then the
same performance can be achieved by using a modified update rule
\begin{equation}
\rho^{(k+1)} = \left\{ \
\begin{array}{ll}
  \tau \rho^{(k)}      & \text{ if } \normsz[\big]{\mb{r}^{(k)}}_2 >
  \xi \mu
  \normsz[\big]{\mb{s}^{(k)}}_2\\[3pt]
  \tau^{-1} \rho^{(k)}  & \text{ if } \normsz[\big]{\mb{s}^{(k)}}_2 >
  \xi^{-1} \mu
  \normsz[\big]{\mb{r}^{(k)}}_2\\[3pt]
  \rho^{(k)}           & \text{ otherwise } \;,
\end{array}
\right.
\label{eq:rhoupdatexi}
\end{equation}
with $\xi = \beta^{-1} \alpha \gamma$ chosen to compensate for the
scaling of the ratio of residuals with the problem scaling.

It is important to emphasise, however, that this issue is not only
relevant to the practitioner considering explicitly scaling an existing
ADMM problem: problem $\tilde{P}$ merely makes explicit the implicit
choices involved in setting up any ADMM problem, and there is no
reason to believe that the often-arbitrary choices made in setting up
the problem correspond to an optimal or even a good choice of scaling
with respect to the convergence of the ADMM iterates subject to update
rule~\eq{rhoupdate}, or subject to update rule~\eq{rhoupdatexi} with
$\xi = 1$.

\subsection{Relative Residuals}
\label{sec:relres}

A simple approach that avoids the need for explicit compensation for
problem scaling when the formulation is modified is to base the
adaptive penalty parameter policy on residuals that represent relative
instead of absolute error\footnote{It is worth noting that similar
  normalisation of error/convergence measures is quite commonly
  applied in other areas of optimization, see e.g.
  \cite[Sec. 1.2]{mittelmann-2003-independent},
  \cite[Sec. 2.1]{wachter-2006-implementation}.}. If the
normalisations required for relative error measures are selected
appropriately\footnote{It is no coincidence that these normalisations
  turn out to be the same as those in the definitions of
  $\epsilon_{\mathrm{pri}}^{(k+1)}$ and
  $\epsilon_{\mathrm{dua}}^{(k+1)}$ in~\cite[Sec
  3.3.1]{boyd-2010-distributed}.}, they will cancel the scaling with
$\beta$ and $\alpha \gamma$, making them invariant to problem
scaling. A reasonable normalisation to make the primal residual
$\mb{r}^{(k+1)} = A \mb{x}^{(k+1)} + B \mb{z}^{(k+1)} - \mb{c}$ a
relative residual is
$$\max\left\{\normsz[\big]{ A \mb{x}^{(k+1)}}_2,
  \normsz[\big]{ B \mb{z}^{(k+1)}}_2, \normsz[\big]{ \mb{c}}_2
\right\} \;,
$$
allowing us to define the relative primal residual
{\small
\begin{equation}
  \mb{r}_{\mathrm{rel}}^{(k+1)} = \frac{A \mb{x}^{(k+1)} + B
    \mb{z}^{(k+1)} - \mb{c}}{\max\left\{\normsz[\big]{ A \mb{x}^{(k+1)}}_2,
    \normsz[\big]{ B \mb{z}^{(k+1)}}_2, \normsz[\big]{ \mb{c}}_2
  \right\}} \;,
\label{eq:nrmprires}
\end{equation}
}
which is invariant to problem scaling since the normalisation factor
has the same scaling as the absolute residual,
{\small
\begin{multline}
\max\left\{\normsz[\big]{\tl{A} \tmb{x}^{(k+1)}}_2,
    \normsz[\big]{\tl{B} \tmb{z}^{(k+1)}}_2, \normsz[\big]{\tmb{c}}_2
  \right\} \\ = \max\left\{\normsz[\big]{\beta \gamma A \gamma^{-1}
      \mb{x}^{(k+1)}}_2, \normsz[\big]{\beta \gamma B \gamma^{-1}
      \mb{z}^{(k+1)}}_2, \normsz[\big]{\beta \mb{c}}_2
  \right\}\\
   = \beta \max\left\{\normsz[\big]{ A \mb{x}^{(k+1)}}_2,
    \normsz[\big]{ B \mb{z}^{(k+1)}}_2, \normsz[\big]{ \mb{c}}_2
  \right\} \;.
\label{eq:rnrmfct}
\end{multline}
}

A suitable normalisation for the dual residual $\mb{s}^{(k+1)} = \rho
A^T B (\mb{z}^{(k+1)} - \mb{z}^{(k)})$ can be obtained
from~\eq{dualres0}.  When $f$ is differentiable and the gradient is
easily computable, a reasonable choice of the normalisation would be
$\max\left\{\normsz[\big]{ \nabla f(\mb{x}^{(k+1)})}_2,
  \normsz[\big]{A^T \mb{y}^{(k+1)}}_2 \right\}$, but since this is
often not the case, we simply use $\normsz[\big]{A^T
  \mb{y}^{(k+1)}}_2$ as the normalisation factor, giving the
relative dual residual
\begin{equation}
  \resizemath{.85\hsize}{
  \mb{s}_{\mathrm{rel}}^{(k+1)} = \frac{\rho
    A^T B (\mb{z}^{(k+1)} \!-\! \mb{z}^{(k)})}{\normsz[\big]{A^T \mb{y}^{(k+1)}}_2
  } = \frac{
    A^T B (\mb{z}^{(k+1)} \!-\! \mb{z}^{(k)})}{\normsz[\big]{A^T \mb{u}^{(k+1)}}_2
  } \;,
  }
\label{eq:nrmduares}
\end{equation}
which is again invariant to problem scaling since the normalisation factor
has the same scaling as the absolute residual,
\vspace{-3mm}
\begin{equation}
\resizemath{.88\hsize}{
\normsz[\big]{\tl{A}^T \tmb{y}^{(k+1)}}_2 = \normsz[\big]{\beta
  \gamma A^T \frac{\alpha}{\beta} \mb{y}^{(k+1)}}_2
  = \alpha \gamma \normsz[\big]{A^T \mb{y}^{(k+1)}}_2 \;.
}
\label{eq:snrmfct}
\end{equation}
Using these definitions, $\tmb{r}_{\mathrm{rel}}^{(k+1)} =
\mb{r}_{\mathrm{rel}}^{(k+1)}$ and $\tmb{s}_{\mathrm{rel}}^{(k+1)} =
\mb{s}_{\mathrm{rel}}^{(k+1)}$; \ie the residuals are invariant to
problem scaling.
The corresponding penalty parameter update policy becomes
\begin{equation}
\rho^{(k+1)} = \left\{ \
\begin{array}{ll}
  \tau \rho^{(k)}      & \text{ if } \normsz[\big]{\mb{r}_{\mathrm{rel}}^{(k)}}_2 >
  \xi \mu
  \normsz[\big]{\mb{s}_{\mathrm{rel}}^{(k)}}_2\\[3pt]
  \tau^{-1} \rho^{(k)}  & \text{ if } \normsz[\big]{\mb{s}_{\mathrm{rel}}^{(k)}}_2 >
  \xi^{-1} \mu
  \normsz[\big]{\mb{r}_{\mathrm{rel}}^{(k)}}_2\\[3pt]
  \rho^{(k)}           & \text{ otherwise } \;,
\end{array}
\right.
\label{eq:rhoupdaterel}
\end{equation}
where the parameter $\xi$ is retained for reasons that will be made
apparent shortly.

The convergence proof~\cite{he-2000-alternating} of the standard
adaptive scheme (\ie~\eq{rhoupdate} with the standard definitions of
the residuals) depends only on bounds on the sequences $\rho^{(k)}$ and
$\eta_k = \sqrt{(\rho^{(k+1)} / \rho^{(k)})^2 - 1}$, neither of which
is affected by the change in the definition of the residuals, so the
convergence results still hold under the modified definitions of the
residuals.

\subsection{Adaptive Multiplier Policy}

The fixed multiplier $\tau$ is a potential weakness of the penalty
update policies~\eq{rhoupdate} and~\eq{rhoupdaterel}. If $\tau$ is
small, then a large number of iterations may be required\footnote{In
  many problems to which ADMM is applied, solving the $\mb{x}$
  update~\eq{admmscaledx} involves solving a large linear system,
  which can be efficiently achieved by pre-computing an LU or Cholesky
  factorization of the system matrix for use in each iteration. Since
  the system matrix depends on $\rho$, it is necessary to re-compute
  the factorization when $\rho$ is updated. (This can be avoided by
  use of an alternative
  factorisation~\cite[Sec. 4.2]{liu-2013-nuclear}, but since this
  method is substantially more computationally expensive in some
  cases, and since a thorough comparison with this alternative is
  beyond the scope of the present paper, it will not be considered
  further here.) Given the computational cost of the factorization, it
  is reasonable to only apply the $\rho$ update at every 10 (for
  example) iterations so that the cost of the factorization can be
  amortised over multiple iterations. This compromise further reduces
  the adaption rate of the adaptive penalty policy.  } to reach an
appropriate $\rho$ value if $\rho^{(0)}$ is poorly chosen (\ie, so that
$\normsz[\big]{\mb{r}^{(0)}}_2 \gg \xi \normsz[\big]{\mb{s}^{(0)}}_2$,
or $\normsz[\big]{\mb{r}^{(0)}}_2 \ll \xi \normsz[\big]{\mb{s}^{(0)}}_2$). On the other hand, if $\tau$ is large, the corrections
to $\rho$ may be too large when $\rho$ is close to the optimal value.

A straightforward solution is to adapt $\tau$ at each iteration
\begin{equation}
\hspace{-2.5mm}\resizemath{.92\hsize}{
\tau^{(k)} = \left\{ \
\begin{array}{ll}
  \!\!\!\sqrt{\xi^{-1} \normsz[\big]{\mb{r}^{(k)}}_2 / \normsz[\big]{\mb{s}^{(k)}}_2} &
   \text{ if }  1 \le
   \sqrt{\xi^{-1} \normsz[\big]{\mb{r}^{(k)}}_2 / \normsz[\big]{\mb{s}^{(k)}}_2}
   < \tau_{\mathrm{max}} \\[3pt]
   \!\!\!\sqrt{\xi \normsz[\big]{\mb{s}^{(k)}}_2 / \normsz[\big]{\mb{r}^{(k)}}_2} &
   \text{ if }  \tau_{\mathrm{max}}^{-1} <
   \sqrt{\xi^{-1} \normsz[\big]{\mb{r}^{(k)}}_2 / \normsz[\big]{\mb{s}^{(k)}}_2}
   < 1 \\[3pt]
   \!\!\!\tau_{\mathrm{max}} & \text{ otherwise } \;,
\end{array}
\right.
}
\label{eq:adapttau}
\end{equation}
where $\tau_{\mathrm{max}}$ provides a bound on $\tau$.  Since $\tau$
is bounded, the convergence results~\cite{he-2000-alternating} still
hold for this extension.

\subsection{Stopping Criteria}
\label{sec:relstopcrit}

The stopping criteria in~\sctn{stopcrit} can be expressed in terms of
the relative residuals $\mb{r}_{\mathrm{rel}}$ and $\mb{s}_{\mathrm{rel}}$
as
\begin{equation}
  \normsz[\big]{\mb{r}_{\mathrm{rel}}^{(k)}}_2 \leq \epsilon_{\mathrm{pri}}^{(k)}
  \;\; \text{ and }
  \normsz[\big]{\mb{s}_{\mathrm{rel}}^{(k)}}_2 \leq \epsilon_{\mathrm{dua}}^{(k)}
\end{equation}
where
\smallmath{
\begin{align}
 \epsilon_{\mathrm{pri}}^{(k)} & \!=\! \sqrt{p} \epsilon_{\mathrm{abs}} /
     \max\left\{\normsz[\big]{A \mb{x}^{(k)}}_2,  \normsz[\big]{B
      \mb{z}^{(k)}}_2, \normsz[\big]{\mb{c}}_2 \right\} \!+\!
 \epsilon_{\mathrm{rel}} \label{eq:eprirel} \\
   \epsilon_{\mathrm{dua}}^{(k)} \!&\! = \sqrt{n} \epsilon_{\mathrm{abs}}  /
      \normsz[\big]{A^T \mb{y}^{(k)}}_2 + \epsilon_{\mathrm{rel}}
       \label{eq:eduarel}  \;.
\end{align}
}
These stopping criteria are invariant to problem scaling when
$\epsilon_{\mathrm{abs}} = 0$.

\subsection{Residual Ratio}
\label{sec:resrat}

While the relative residuals proposed in~\sctn{relres}
address the absence of scaling invariance in the adaptive
penalty parameter strategy based on residual balancing, there is
another even more serious deficiency that is not so easily remedied.
As discussed in~\sctn{adaptrho}, the target ratio of unity is
motivated by representing the distance from convergence as $\|
\mb{r}^{(k)} \|_2^2 + \| \mb{s}^{(k)}\|_2^2$, but this greatly
simplifies the true picture.

The ADMM convergence proof in~\cite{boyd-2010-distributed} (see
Sec. 3.3.1 and Appendix A) provides some insight into the relationship
between the distance from convergence and the residuals, in the form of
the inequality
\begin{align}
f(\mb{x}^{(k)}) + g(\mb{z}^{(k)}) - p^* \leq \; & - (\mb{y}^{(k)})^T
\mb{r}^{(k)}  \nonumber \\ & + (\mb{x}^{(k)} - \mb{x}^*)^T \mb{s}^{(k)} \;,
\end{align}
which implies the looser inequality
\begin{align}
f(\mb{x}^{(k)}) + g(\mb{z}^{(k)}) - p^* \leq \; &  \normsz[\big]{ \mb{y}^{(k)} }
 \normsz[\big]{ \mb{r}^{(k)} } + \nonumber \\ & \normsz[\big]{\mb{x}^{(k)} -
   \mb{x}^*} \normsz[\big]{\mb{s}^{(k)}}
\label{eq:convineq}
\end{align}
in terms of the norms of the relevant vectors. Applying the original
argument that led to unity as the appropriate target ratio to this
inequality implies that the appropriate ratio is, in fact,
approximately
$\normsz[\big]{ \mb{y}^{(k)} } / \normsz[\big]{\mb{x}^{(k)} -
  \mb{x}^*}$.
This would explain why some authors have found the original residual
balancing strategy of~\eq{rhoupdate} to be effective~\cite{hansson-2012-subspace, liu-2013-nuclear, vu-2013-fantope, iordache-2014-collaborative, weller-2014-phase, wohlberg-2014-efficient} and others have not~\cite[Sec. 2.4]{ramdas-2015-fast}: the method succeeds when this ratio happens to be relatively close to unity, and fails when it is not.

Unfortunately, since $\mb{x}^*$ is unknown while solving the problem, there is no obvious way to estimate this ratio, and we are left with the rather unsatisfactory solution of accepting $\xi$ in~\eq{rhoupdaterel} as a user-selected parameter of the method. Since this approach essentially replaces one user parameter, $\rho$, with another, $\xi$, it is not clear that the residual balancing strategy has any real value as a parameter selection technique. One might argue that, since the residual balancing method has been found to be satisfactory in a variety of applications, it must often be the case that $\xi = 1$ is not too far from the optimal setting, and that $\xi$ may be a more stable parameterisation than $\rho$, but further study is necessary before any reliable conclusions can be drawn.

Since~\eq{rhoupdaterel} retains $\xi$, which can be used to compensate for explicit problem scaling as discussed in~\sctn{adaptpenparam}, it is reasonable to ask whether there is any real benefit to using~\eq{rhoupdaterel} based on the relative residuals, \ie, since we have an unknown $\xi$ in both cases, what is the advantage of one scaling of this unknown quantity in comparison with another. Two arguments can be made in favour of the use of relative residuals as in~\eq{rhoupdaterel}:
\begin{itemize}
\item Ignoring the question of determining a good choice of $\xi$, once one has been found, \eq{rhoupdaterel} is invariant to problem scaling, while~\eq{rhoupdate} is not.
\item Since~\eq{rhoupdaterel} is invariant to problem scaling, one might expect that the $\xi$ for this update rule is more stable than the $\xi$ for~\eq{rhoupdaterel}, in the sense that it varies across a smaller numerical range for different problem. (This important question is not explored in the experimental results presented here.)
\end{itemize}

It should also be noted that the unknown scalings of the residuals in~\eq{convineq} imply that neither the absolute nor relative stopping tolerances in~\sctn{relstopcrit} can be viewed as providing an actual bound on the solution optimality, either in an absolute or a relative sense (e.g. a relative stopping criterion $\epsilon_{\mathrm{rel}} = 10^{-3}$ does \emph{not} imply that the final iterate is within $10^{-3}$ relative distance to the optimal solution).

\section{BPDN}
\label{sec:bpdn}

To illustrate these issues, we will focus on Basis Pursuit DeNoising
(BPDN)~\cite{chen-1998-atomic},
\begin{equation}
\argmin_{\mb{x}} \frac{1}{2} \norm{D \mb{x} - \bvsigma}_2^2 + \lambda
\norm{\mb{x}}_1 \;,
\label{eq:bpdn}
\end{equation}
a standard problem in computing sparse representations corresponding
to~\eq{admmprob} with
\begin{gather}
f(\mb{x}) = \frac{1}{2} \norm{D \mb{x} - \bvsigma}_2^2 \;\;\;\;
g(\mb{z}) = \lambda \norm{\mb{z}}_1 \nonumber \\
A = I \;\;\;\; B = -I \;\;\;\; \mb{c} = 0 \;.
\end{gather}
Solving via ADMM, we have problem $P$
\begin{equation}
\argmin_{\mb{x}} \frac{1}{2} \norm{D \mb{x} - \bvsigma}_2^2 + \lambda
\norm{\mb{z}}_1 \text{ s.t. } \mb{x} = \mb{z}
\end{equation}
with Lagrangian
\begin{equation}
  L(\mb{x}, \mb{z}, \mb{y}) = \frac{1}{2} \norm{D \mb{x} -
    \bvsigma}_2^2 + \lambda \norm{\mb{z}}_1 + \mb{y}^T (\mb{x} - \mb{z}) \;.
\end{equation}

We also consider Convolutional BPDN (CBPDN), a variant of BPDN
constructed by replacing the linear combination of a set of dictionary
vectors by the sum of a set of convolutions with dictionary
filters~\cite[Sec. II]{wohlberg-2016-efficient}
\begin{equation}
\argmin_{\{\mb{x}_m\}} \frac{1}{2} \normsz[\Big]{\sum_m \mb{d}_m \ast \mb{x}_m
- \bvsigma}_2^2 + \lambda \sum_m \norm{\mb{x}_m}_1 \; ,
\label{eq:convbpdn}
\end{equation}
where $\{\mb{d}_m\}$ is a set of $M$ dictionary \emph{filters}, $\ast$
denotes convolution, and $\{\mb{x}_m\}$ is a set of coefficient
maps. Algebraically, this variant is a special case of standard BPDN, so
that the same scaling properties apply, but since the dictionaries in
this form are very highly overcomplete (the overcompleteness factor is
equal to the number of filters $M$), one may expect that this variant
might exhibit at least somewhat different behaviour in practice. A
further difference is that the $\{\mb{x}_m\}$ can be efficiently
computed without any factorisation of system
matrices~\cite{wohlberg-2014-efficient}, so in this case the penalty
update policy is applied at every iteration instead of at every 10
iterations.

\section{Results}
\label{sec:rslt}

In this section the issues discussed above are illustrated via a number of computational experiments. Many of these experiments compare the effect of different penalty parameter selection methods on the number of iterations required to reach the stopping criteria. With respect to these experiments, it must be emphasised that:
\begin{itemize}
\item Since the relationship between the stopping criteria and the actual solution suboptimality is unknown (see~\sctn{resrat}), reaching the stopping criteria faster does \emph{not} imply faster convergence.
\item These experiments all use relative stopping thresholds (\ie $\epsilon_{\mathrm{abs}} = 0$), which could be considered to confer an advantage on the relative residual balancing policy since it balances the residuals in a way that that is favourable to satisfying the relative stopping thresholds\footnote{Since the stopping criteria require that both residuals are below the same threshold, they will be satisfied more quickly if they are roughly equal than if one is much larger than the other, all else being equal.}. Note, however, that the original goal of invariance to problem scaling cannot be achieved if $\epsilon_{\mathrm{abs}} \neq 0$.
\end{itemize}

\subsection{BPDN with Random Dictionary}

\begin{figure}[htbp]
  \small
  \hspace{-5mm}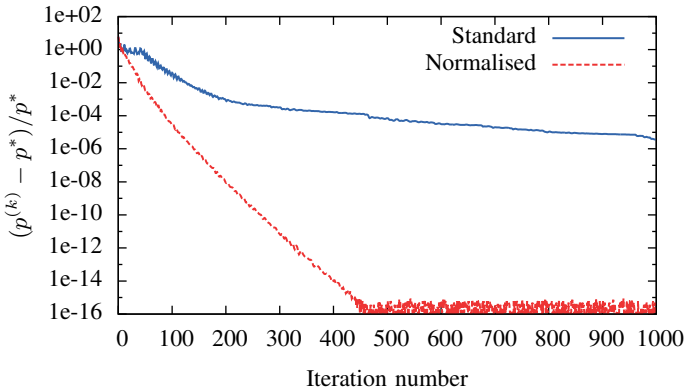
  \vspace{-3mm}
  \caption{A comparison of functional value evolution for the same
    problem with adaptive $\rho$ based on standard and normalised residuals.}
  \label{fig:exp1fnval}
\end{figure}

\begin{figure}[htbp]
  \small
  \hspace{-5mm}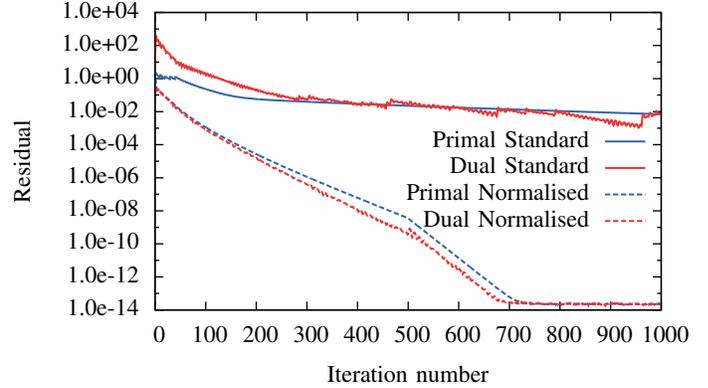
  \vspace{-3mm}
  \caption{A comparison of primal and dual residual evolution for the
    same problem with adaptive $\rho$ based on standard and normalised
    residuals. For a meaningful comparison, the residuals are divided
    by their respective values of $\epsilon_{\mathrm{pri}}$ or
    $\epsilon_{\mathrm{dua}}$.
  }
  \label{fig:exp1priduares}
\end{figure}

\begin{figure}[htbp]
  \small
  \hspace{-5mm}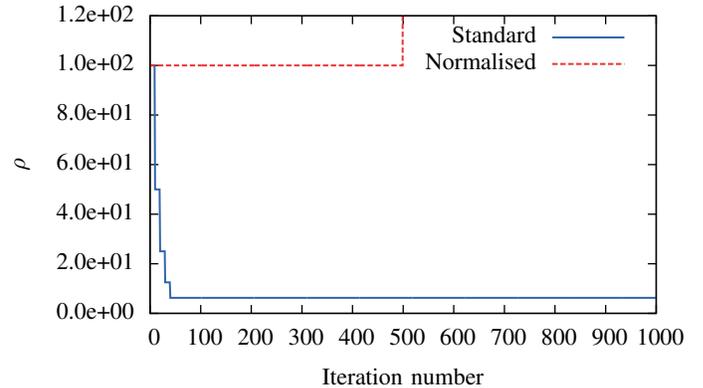
  \vspace{-3mm}
  \caption{A comparison of selected $\rho$ values for the same
    problem with adaptive $\rho$ based on standard and normalised residuals.}
  \label{fig:exp1rho}
\end{figure}

The first experiment involves sparse coefficient recovery on a random
dictionary without normalisation. A dictionary $D \in \mbb{R}^{512
  \times 4096}$ was generated with unit standard deviation
i.i.d. entries with a Gaussian distribution, a corresponding reference
coefficient vector $\mb{x}_0$ was constructed by assigning random
values to 64 randomly selected coefficients, the remainder of which
were zero, and a test signal was constructed by adding Gaussian white
noise of standard deviation 0.5 to the product of $D$ and
$\mb{x}_0$. The experiment involves using BPDN with $\lambda = 40$
(selected for good support identification), $\xi = 1$, $\epsilon_{\mathrm{abs}} =
0$, and $\epsilon_{\mathrm{rel}} = 10^{-4}$ to attempt to recover
$\mb{x}_0$ from the signal, comparing performance with both standard
and normalised residuals. It is clear from
\figs{exp1fnval}--\fign{exp1rho} that the adaptive $\rho$ policy gives
very substantially better performance with normalised residuals than
with the standard definition. The desired stopping tolerance is
reached within 160 iterations when using normalised residuals, but has
still not been attained when the maximum iteration limit of 1000 is
reached in the case of standard residuals. The performance difference
is even greater if random dictionary $D$ is generated with standard
deviation greater than unity.

\subsection{BPDN with Learned Dictionary}

\begin{figure*}[htbp]
  \centering \small
  \begin{tabular}{cc}
    \hspace{-10mm}\subfloat[\label{fig:bpdn_exp12_rhoplot1_std}
    \protect\rule{0pt}{1.5em}
    Standard residuals]
               {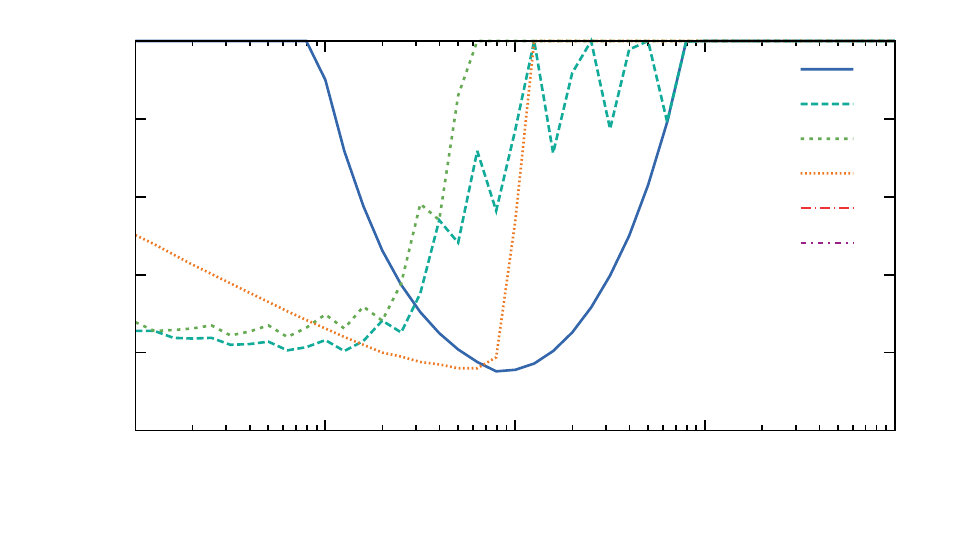} &
   \hspace{-9mm}\subfloat[\label{fig:bpdn_exp12_rhoplot1_nrm}
   \protect\rule{0pt}{1.5em}
    Normalised residuals]
               {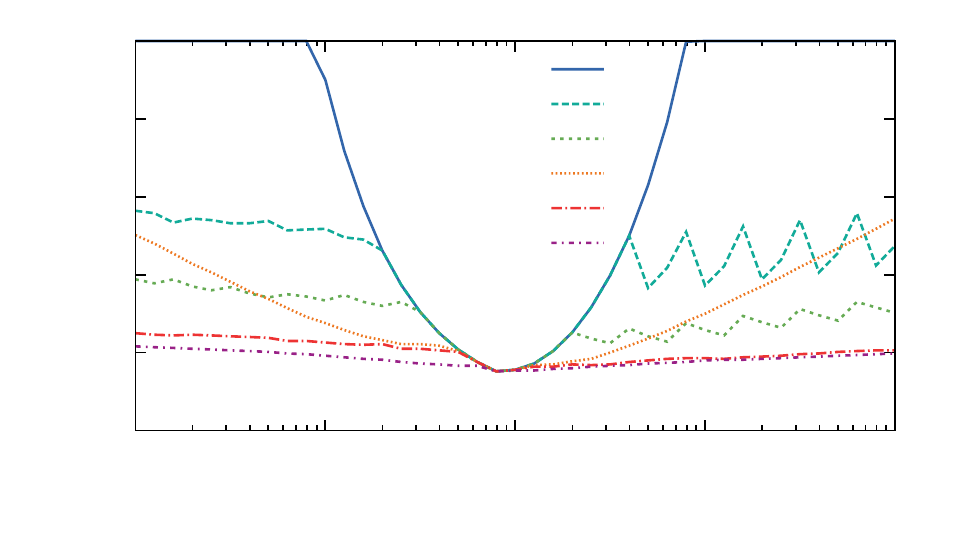}
  \end{tabular}
  \caption{Variation with $\rho^{(0)}$ of number of iterations
    required to reach a relative stopping tolerance of
    $\epsilon_{\mathrm{rel}} = 10^{-3}$ for different variants of the
    adaptive $\rho$ policy, and for standard and normalised residuals,
    in a BPDN problem with $D \in \mbb{R}^{64 \times 128}$ and
    $\lambda = 10^{-2}$. The variant labels are ``Fixed'', indicating
    that $\rho$ is fixed at $\rho^{(0)}$ and is not adapted, of the
    form $\mu /\tau$, or of the form $\mu$/Auto, which indicates that
    $\tau$ is adapted as in~\eq{adapttau}, with $\tau_{\mathrm{max}} =
    100$ and $\xi = 1$.}
  \label{fig:bpdn_exp12_rhoplot1}
\end{figure*}

\begin{figure*}[htbp]
  \centering \small
  \begin{tabular}{cc}
    \hspace{-10mm}\subfloat[\label{fig:bpdn_exp12_64_std}
    \protect\rule{0pt}{1.5em}
    Standard residuals]
               {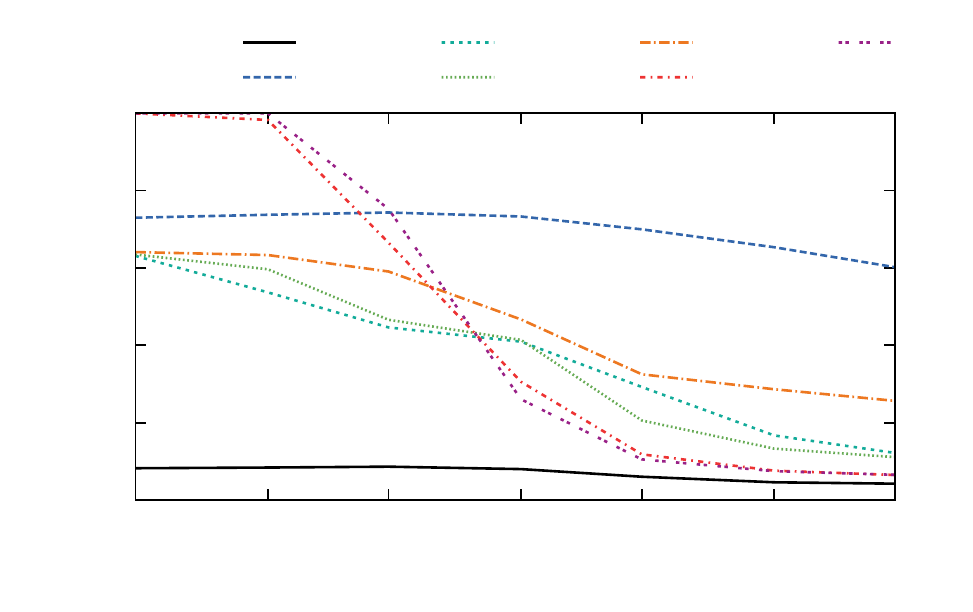} &
   \hspace{-9mm}\subfloat[\label{fig:bpdn_exp12_64_nrm}
   \protect\rule{0pt}{1.5em}
    Normalised residuals]
               {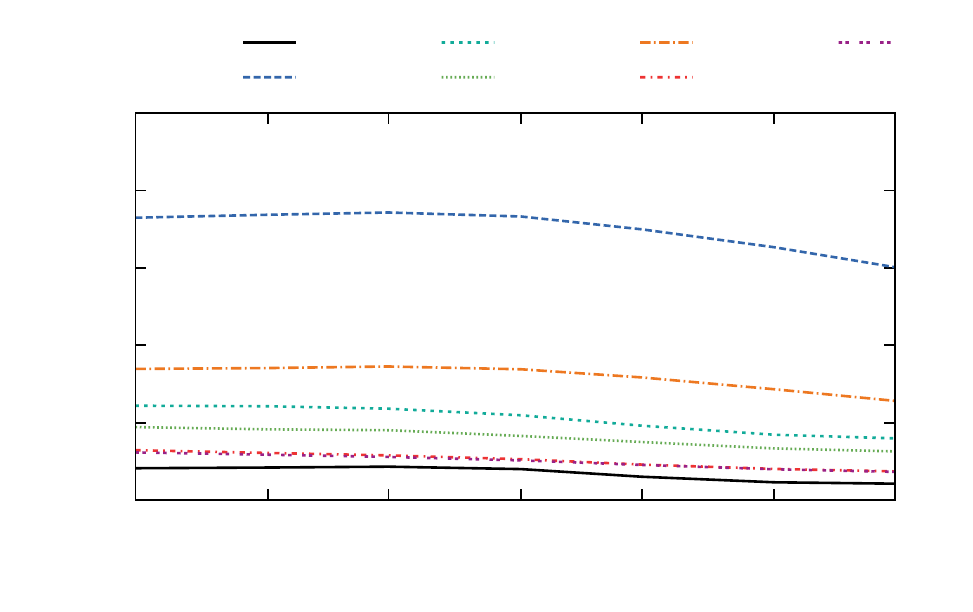}
  \end{tabular}
  \caption{Mean number of iterations (averaged over all values of
    $\rho^{(0)}$) required to reach a relative stopping tolerance of
    $\epsilon_{\mathrm{rel}} = 10^{-3}$ for different variants of the
    adaptive $\rho$ policy, and for standard and normalised residuals,
    in a BPDN problem with $D \in \mbb{R}^{64 \times 64}$ and varying
    $\lambda$. The variant labels are ``Fixed'', indicating that
    $\rho$ is fixed at $\rho^{(0)}$ and is not adapted, of the form
    $\mu /\tau$, or of the form $\mu$/Auto, which indicates that
    $\tau$ is adapted as in~\eq{adapttau}, with $\tau_{\mathrm{max}} =
    100$ and $\xi = 1$. ``Fixed (min)'' denotes the minimum number of
    iterations (i.e. not the mean) obtained via the best fixed choice
    of $\rho^{(0)}$ at each value of $\lambda$.}
  \label{fig:bpdn_exp12_64}
\end{figure*}

\begin{figure*}[htbp]
  \centering \small
  \begin{tabular}{cc}
    \hspace{-10mm}\subfloat[\label{fig:bpdn_exp12_96_std}
    \protect\rule{0pt}{1.5em}
    Standard residuals]
               {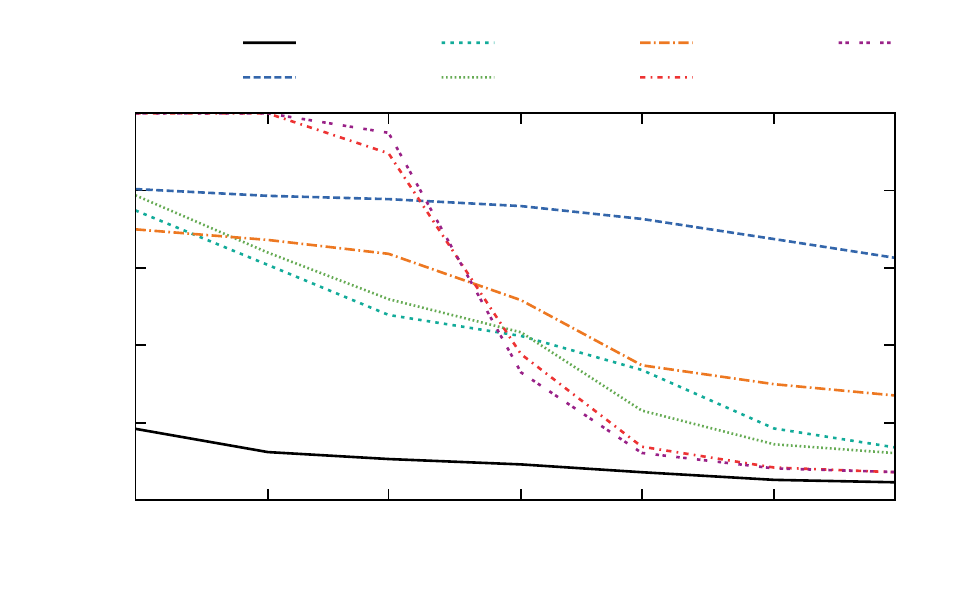} &
   \hspace{-9mm}\subfloat[\label{fig:bpdn_exp12_96_nrm}
   \protect\rule{0pt}{1.5em}
    Normalised residuals]
               {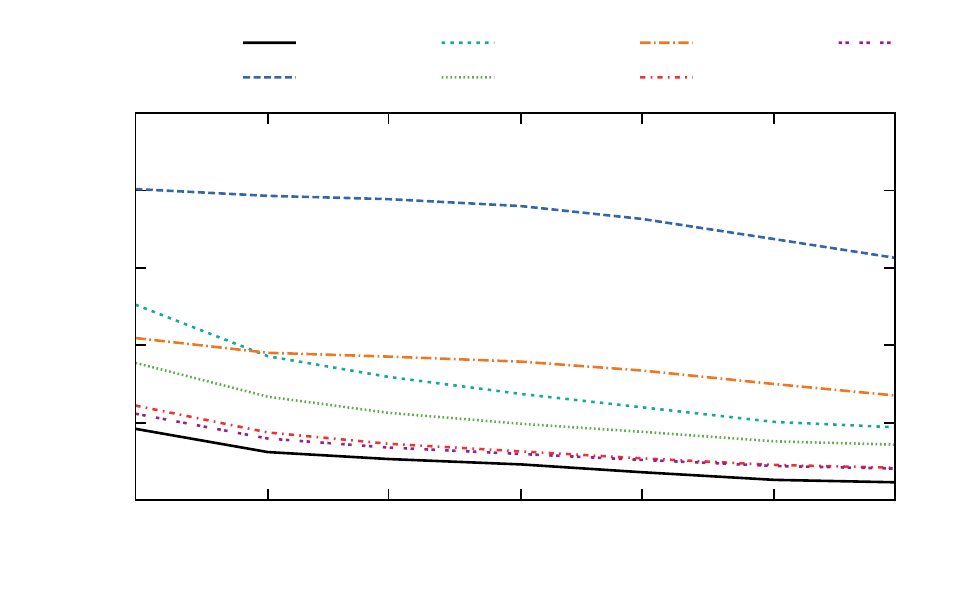}
  \end{tabular}
  \caption{Mean number of iterations (averaged over all values of
    $\rho^{(0)}$) required to reach a relative stopping tolerance of
    $\epsilon_{\mathrm{rel}} = 10^{-3}$ for different variants of the
    adaptive $\rho$ policy, and for standard and normalised residuals,
    in a BPDN problem with $D \in \mbb{R}^{64 \times 96}$ and varying
    $\lambda$. The variant labels are ``Fixed'', indicating that
    $\rho$ is fixed at $\rho^{(0)}$ and is not adapted, of the form
    $\mu /\tau$, or of the form $\mu$/Auto, which indicates that
    $\tau$ is adapted as in~\eq{adapttau}, with $\tau_{\mathrm{max}} =
    100$ and $\xi = 1$. ``Fixed (min)'' denotes the minimum number of
    iterations (i.e. not the mean) obtained via the best fixed choice
    of $\rho^{(0)}$ at each value of $\lambda$.}
  \label{fig:bpdn_exp12_96}
\end{figure*}

\begin{figure*}[htbp]
  \centering \small
  \begin{tabular}{cc}
    \hspace{-10mm}\subfloat[\label{fig:bpdn_exp12_128_std}
    \protect\rule{0pt}{1.5em}
    Standard residuals]
               {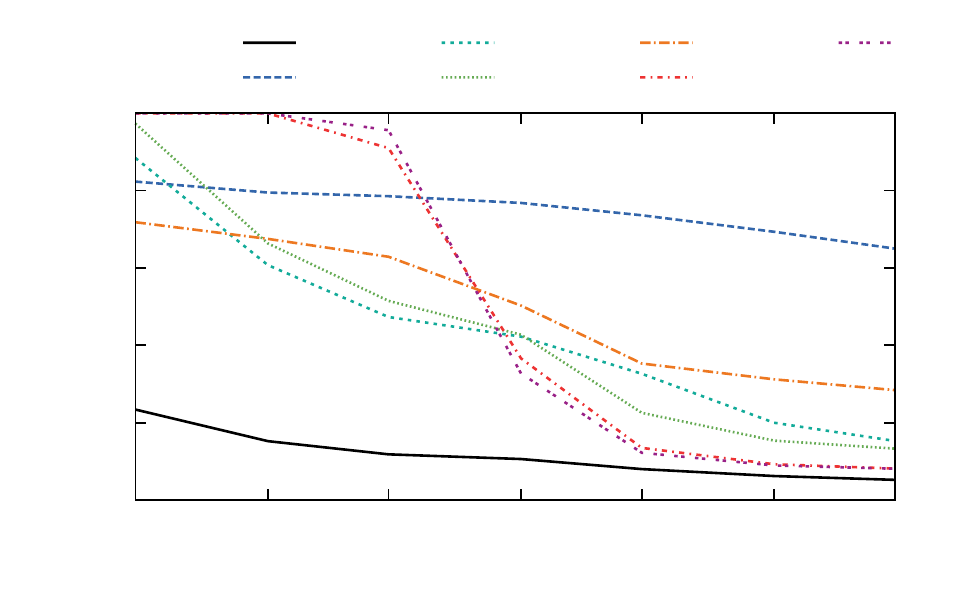} &
   \hspace{-9mm}\subfloat[\label{fig:bpdn_exp12_128_nrm}
   \protect\rule{0pt}{1.5em}
    Normalised residuals]
               {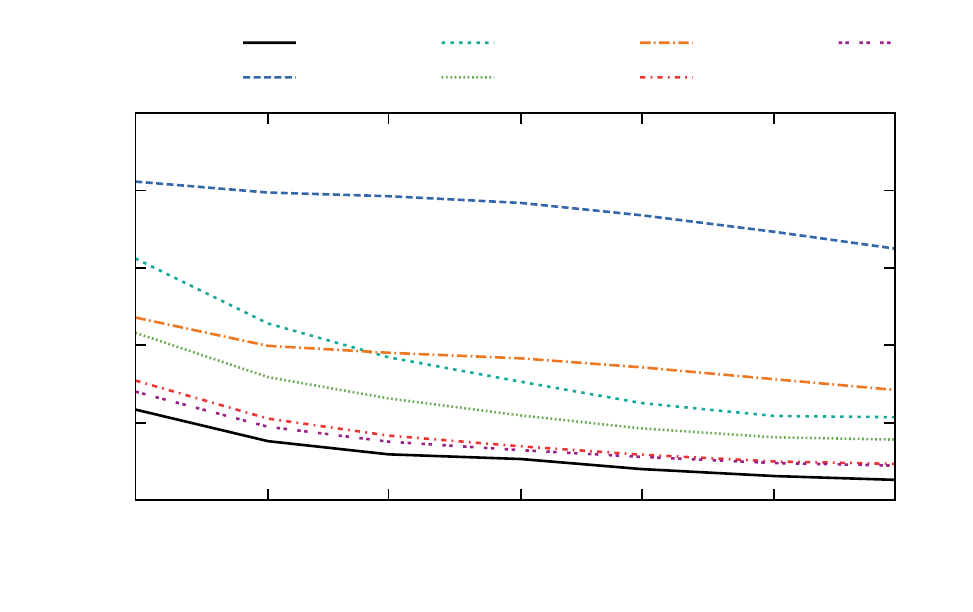}
  \end{tabular}
  \caption{Mean number of iterations (averaged over all values of
    $\rho^{(0)}$) required to reach a relative stopping tolerance of
    $\epsilon_{\mathrm{rel}} = 10^{-3}$ for different variants of the
    adaptive $\rho$ policy, and for standard and normalised residuals,
    in a BPDN problem with $D \in \mbb{R}^{64 \times 128}$ and varying
    $\lambda$. The variant labels are ``Fixed'', indicating that
    $\rho$ is fixed at $\rho^{(0)}$ and is not adapted, of the form
    $\mu /\tau$, or of the form $\mu$/Auto, which indicates that
    $\tau$ is adapted as in~\eq{adapttau}, with $\tau_{\mathrm{max}} =
    100$ and $\xi = 1$. ``Fixed (min)'' denotes the minimum number of
    iterations (i.e. not the mean) obtained via the best fixed choice
    of $\rho^{(0)}$ at each value of $\lambda$.}
  \label{fig:bpdn_exp12_128}
\end{figure*}

The second set of experiments compares the performance of a fixed
$\rho$ and various adaptive $\rho$ parameter choices, using standard
and normalised residuals, for a Multiple Measurement Vector (MMV) BPDN
problem. Dictionaries $D \in \mbb{R}^{64 \times 64}$, $D \in
\mbb{R}^{64 \times 96}$, and $D \in \mbb{R}^{64 \times 128}$ were
learned on a large training set of $8 \times 8$ image patches, and the
test data consisted of 32558 zero-mean $8 \times 8$ image patches
represented as a matrix $S \in \mbb{R}^{64 \times 32258}$. The number
of iterations required to attain a relative stopping tolerance of
$\epsilon_{\mathrm{rel}} = 10^{-3}$ for $D \in \mbb{R}^{64 \times
  128}$ and $\lambda = 10^{-2}$ is
compared in~\fig{bpdn_exp12_rhoplot1}.
The following observations can be made
with respect to the ability of the different methods to reduce the
dependence of the number of iterations on the initial choice
$\rho^{(0)}$:
\begin{itemize}
\item The best choice of fixed $\rho$ gives similar performance to
  the best adaptive strategy, but performance falloff is quite
  rapid as $\rho$ is changed away from the optimum. Given the
  absence of techniques for identifying the optimum $\rho$ \emph{a
    priori} for most problems, it is clear that the adaptive
  strategy can play a valuable role in reducing computation time.
\item When using normalised residuals, there is an overall improvement
 with smaller $\mu$. In particular, it appears that, at least for the
 BPDN problem, the standard choice of $\mu = 10$ is too coarse, and
 benefit can be obtained from finer control of the residual ratio,
\item When using standard residuals, the converse is true, performance
  decreasing with smaller $\mu$. This should not be surprising given
  the previously identified theoretical problems regarding the use of
  standard residuals in~\eq{rhoupdate}: the errors in the residual
  ratio that are masked by setting $\mu = 10$ become increasingly
  apparent as $\mu$ is reduced in an attempt at exerting finer control
  over the residual ratio. In this case the performance of the
  adaptive $\tau$ methods based on~\eq{adapttau} is particularly poor
  because the adaptive $\tau$ allows $\rho$ to be more rapidly
  adjusted to the incorrect value based on the incorrect residual
  ratios.
\item The best overall performance is provided by the two automatic
  $\tau$ methods based on~\eq{adapttau} with normalised residuals.
\end{itemize}

Comparisons of the different strategies over a wide range of $\lambda$
values and three different dictionary sizes are presented
in~\figs{bpdn_exp12_64}--\fign{bpdn_exp12_128}. The mean number of
iterations for all $\rho$ values is plotted against $\lambda$, and
also compared with the minimum number of iterations obtained for the
best fixed choice of $\rho$. The most important observations to be
made are:
\begin{itemize}
\item The standard residuals give similar performance to the
  normalised residuals for the larger values of $\lambda$ since in
  this regime the normalisation quantities turn out to be close to
  unity.
\item At smaller values of $\lambda$, the normalised residuals give
  much better performance.
\item Considered over the entire range of $\lambda$ values, the
  normalised residuals all give better performance than their
  un-normalised counterparts.
 \item Of the methods using normalised residuals, the adaptive $\tau$
   methods based on~\eq{adapttau} gives substantially better
   performance than the standard methods.
\end{itemize}

\subsection{Convolutional BPDN Problem}
\label{sec:cbpdnrslt}

\begin{figure*}[htbp]
  \centering \small
  \begin{tabular}{cc}
    \hspace{-10mm}\subfloat[\label{fig:cbpdn_exp04_64_std}
    \protect\rule{0pt}{1.5em}
    Standard residuals]
               {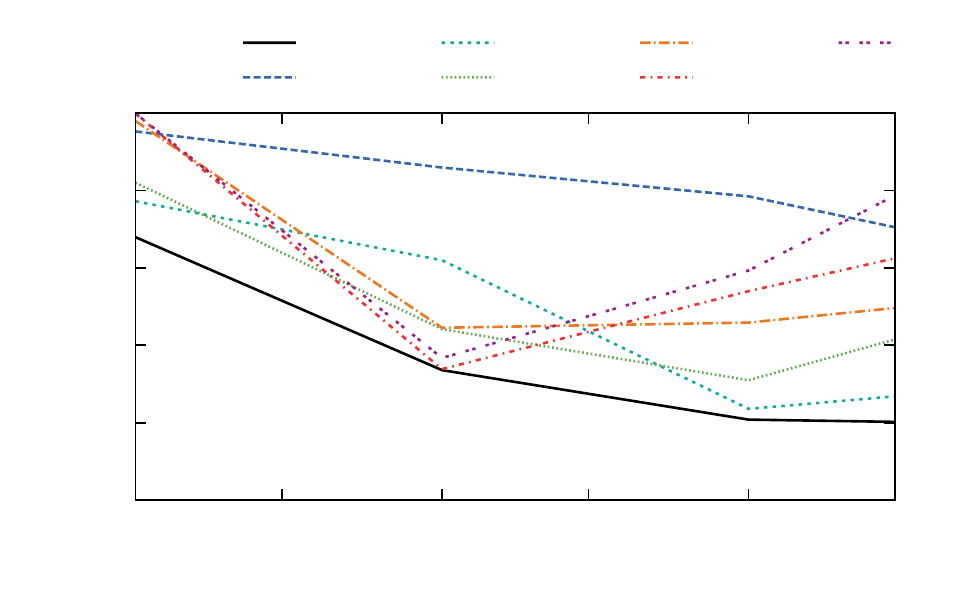} &
   \hspace{-9mm}\subfloat[\label{fig:cbpdn_exp04_64_nrm}
   \protect\rule{0pt}{1.5em}
    Normalised residuals]
               {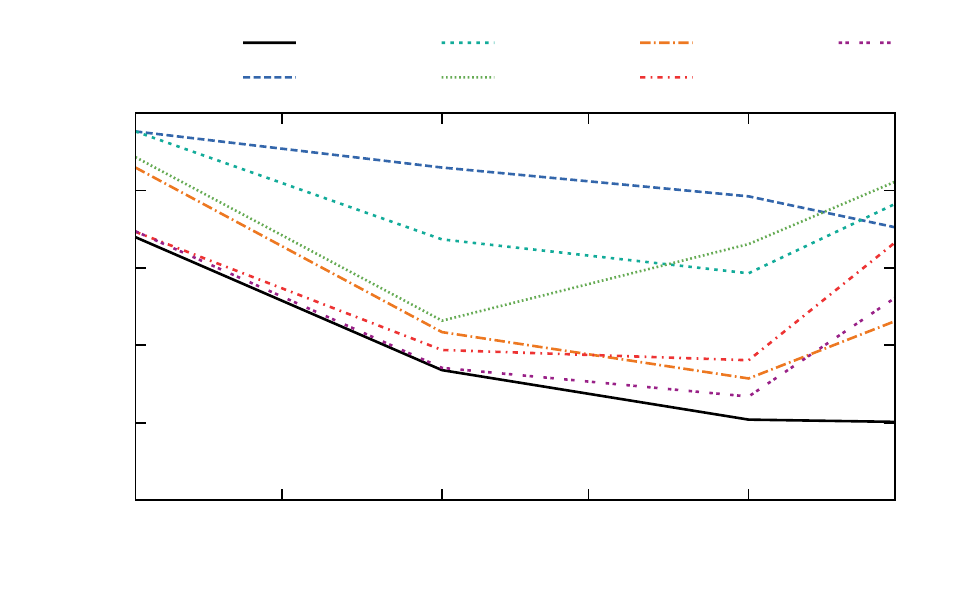}
  \end{tabular}
  \caption{Mean number of iterations (averaged over all values of
    $\rho^{(0)}$) required to reach a relative stopping tolerance of
    $\epsilon_{\mathrm{rel}} = 10^{-3}$ for different variants of the
    adaptive $\rho$ policy, and for standard and normalised residuals,
    in a CBPDN problem with a $8 \times 8 \times 64$ dictionary and
    varying $\lambda$. The variant labels are ``Fixed'', indicating
    that $\rho$ is fixed at $\rho^{(0)}$ and is not adapted, of the
    form $\mu /\tau$, or of the form $\mu$/Auto, which indicates that
    $\tau$ is adapted as in~\eq{adapttau}, with $\tau_{\mathrm{max}} =
    100$ and $\xi = 1$. ``Fixed (min)'' denotes the minimum number of
    iterations (i.e. not the mean) obtained via the best fixed choice
    of $\rho^{(0)}$ at each value of $\lambda$.}
  \label{fig:cbpdn_exp04_64}
\end{figure*}

The penalty update strategies were also compared in application to a
Convolutional BPDN problem consisting of jointly computing the
representations of two $256 \times 256$ pixel images\footnote{As is
  common practice in convolutional sparse representations, the
  representation was computed after a highpass filtering
  pre-processing step, consisting in this case of application of a
  lowpass filter, equivalent to solving the problem $\argmin_{\mb{x}}
  \frac{1}{2} \norm{\mb{x} - \bvsigma}_2^2 + \lambda_L \norm{\nabla
    \mb{x}}_2^2$ with $\lambda_L = 5.0$, and then subtracting the
  lowpass filtered images from the corresponding original images.}
(the well-known ``Lena'' and ``Barbara'' images), with a dictionary
consisting of 64 filters of size $8 \times 8$ samples and for a range
of $\lambda$ and $\rho^{(0)}$ values. It can be seen
from~\fig{cbpdn_exp04_64} that the normalised residuals give good
performance for $\lambda \leq 0.1$, but for larger values of $\lambda$
neither standard nor normalised residuals provide performance close to
that of the best fixed $\rho$.
\begin{figure}[htbp]
  \small
  \hspace{-5mm}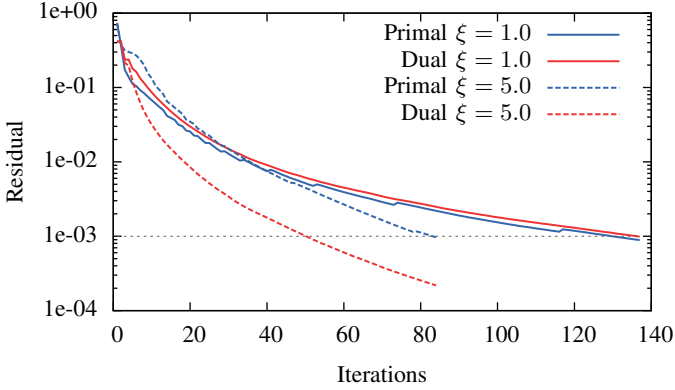
  \vspace{-2mm}
  \caption{Evolution of primal and dual residuals for two different
    choices of $\xi$ in a CBPDN problem with an $8 \times 8 \times 32$
    dictionary, $\lambda = 0.3$, and $\rho^{(0)} = 251$. The $\rho$
    update policy was as in~\eq{rhoupdaterel}, with normalised
    residuals, $\mu = 1.2$, and with adaptive $\tau$ as
    in~\eq{adapttau}, with $\tau_{\mathrm{max}} = 100$.}
  \label{fig:cbpdnexp09rsxmpl}
\end{figure}
This is an indication that $\xi = 1$ is not a suitable choice in this case, for which $\xi = 5$ gives better performance, as illustrated in~\fig{cbpdnexp09rsxmpl}.

\begin{figure}[htbp]
  \small
  \hspace{-1mm}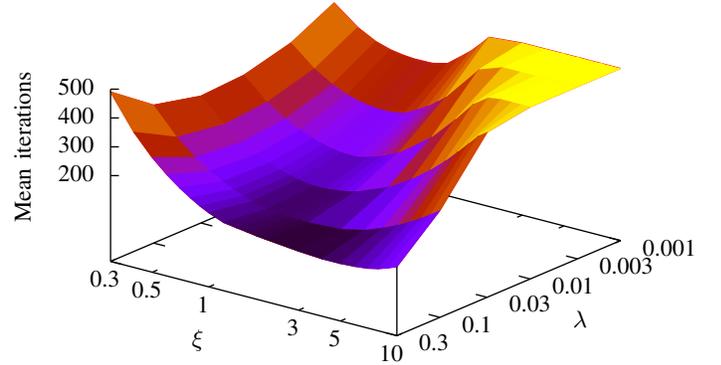
  \vspace{-6mm}
  \caption{Mean number of iterations, averaged over all values of $\rho^{(0)}$,
    against $\lambda$ and $\xi$ for a CBPDN problem with
    an $8 \times 8 \times 64$ dictionary.  The $\rho$
    update policy was as in~\eq{rhoupdaterel}, with normalised
    residuals, $\mu = 1.2$, and with adaptive $\tau$ as
    in~\eq{adapttau}, with $\tau_{\mathrm{max}} = 1000$.}
  \label{fig:cbpdnexp02itlmdxi}
\end{figure}

\begin{figure}[htbp]
  \small
  \hspace{-1mm}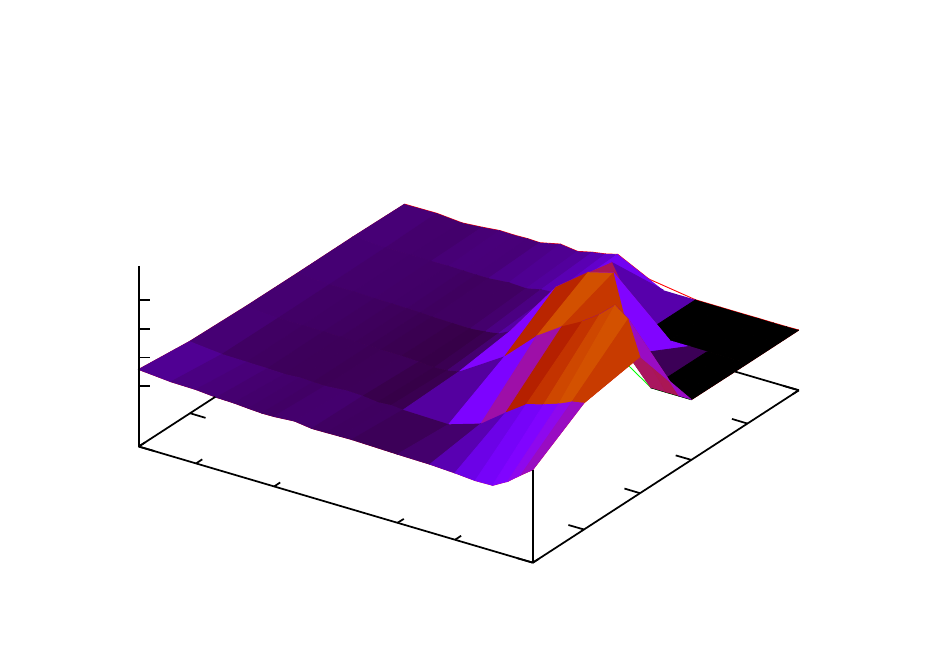
  \vspace{-6mm}
  \caption{Standard deviation of number of iterations with respect to
     $\rho^{(0)}$ in~\fig{cbpdnexp02itlmdxi}. Note that the variation with
    respect to $\rho$ is small where the mean number of iterations is
    small. The standard deviation is zero for small $\lambda$ and
  large $\xi$ because the number of iterations is clipped to 500 by
  the maximum iteration limit in this region.}
  \label{fig:cbpdnexp02sditlmdxi}
\end{figure}

The effect of varying $\xi$ was investigated by running a large number
of computational experiments for the CBPDN problem, with a $8 \times 8
\times 64$ dictionary and for different values of $\lambda$ (6
approximately logarithmically spaced values in the range $1 \times
10^{-3}$ to $0.3$), $\rho$ (51 logarithmically spaced values in the
range $10^{-1}\lambda$ to $10^{4}\lambda$), and $\xi$ (21 values in
the range 0.3 to 10.0).  The mean and standard deviation over
$\rho^{(0)}$ of the number of iterations required to reach stopping
tolerance $\epsilon_{\mathrm{abs}} = 0, \epsilon_{\mathrm{rel}} =
10^{-3}$ are displayed in~\fig{cbpdnexp02itlmdxi}
and~\fign{cbpdnexp02sditlmdxi} respectively. It can be observed that
that the value of $\xi$ giving the minimum number of iterations varies
with $\lambda$, and that considering the mean over $\rho^{(0)}$ of the
number of iterations is a reasonable criterion since the variation
with $\rho^{(0)}$ is small when $\xi$ is well chosen.

\begin{figure*}[htbp]
  \centering \small
  \begin{tabular}{cc}
    \hspace{-10mm}\subfloat[\label{fig:cbpdnexp09xi}\protect\rule{0pt}{1.5em}
    Function fit to best values of $\xi$ for different $\lambda$.
    ]
               {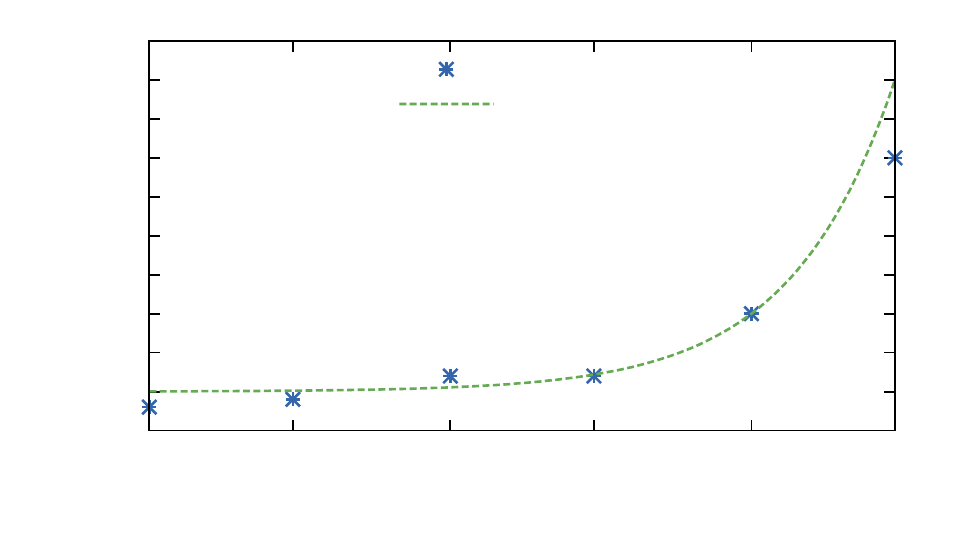} &
   \hspace{-9mm}\subfloat[\label{fig:cbpdnexp09it}\protect\rule{0pt}{1.5em}
    Mean iterations for different values of $\lambda$.
    ]
               {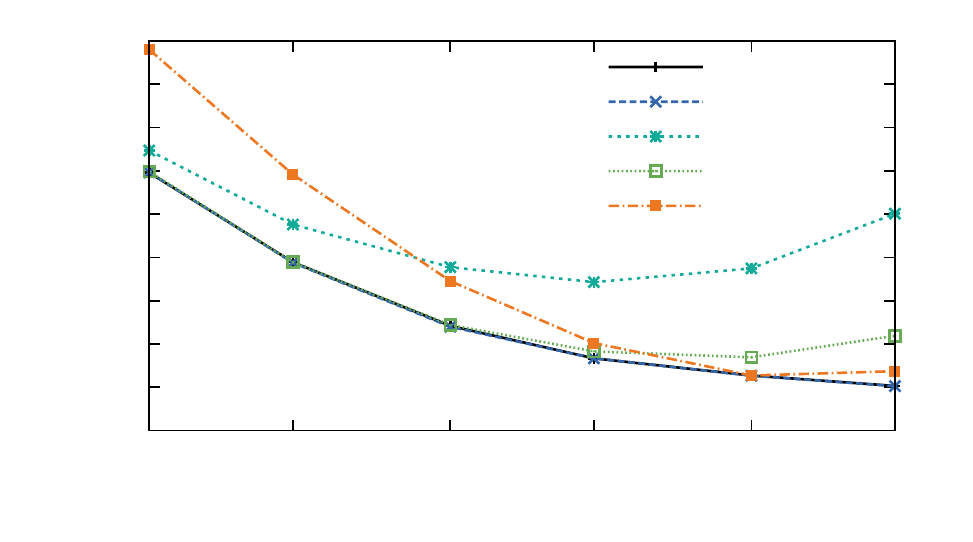}
  \end{tabular}
  \caption{(a) shows the good fit of function $f(\lambda) = 1 +
    18.3^{\log_{10}(\lambda) + 1}$ to the values of $\xi$ that
    minimise the mean (over all values of $\rho^{(0)}$) number of
    required iterations for different values of $\lambda$, determined
    by running a large number of simulations for different values of
    $\xi$, $\rho$, and $\lambda$. (b) shows the variation with
    $\lambda$ of the mean (over all values of $\rho^{(0)}$) number of
    iterations for the best choice of $\xi$ as in (a), for $\xi$
    chosen according to the function $f(\lambda)$, and for three fixed
    choices of $\xi$. All simulations were for a CBPDN problem with a
    $8 \times 8 \times 64$ dictionary and $\epsilon_{\mathrm{abs}} =
    0, \epsilon_{\mathrm{rel}} = 10^{-3}$.  The $\rho$ update policy
    was as in~\eq{rhoupdaterel}, with normalised residuals, $\mu =
    1.2$, and with adaptive $\tau$ as in~\eq{adapttau}, with
    $\tau_{\mathrm{max}} = 1000$.}
  \label{fig:cbpdnexp09}
\end{figure*}

Since the best $\xi$ varies with $\lambda$, it is reasonable to ask,
in the absence of any theory to guide the choice, whether there is a
reliable way of making a good choice of $\xi$. By examining the data
for the experiments used to generate~\fig{cbpdnexp02itlmdxi}
and~\fign{cbpdnexp02sditlmdxi}, as well as for corresponding
experiments with other dictionaries with 32, 96, and 128 filters of
size $8 \times 8$, it was determined that the function $f(\lambda) = 1 +
a^{\log_{10}(\lambda) + 1}$ with $a = 18.3$ provides a reasonable fit to
the best choice of $\xi$ for each $\lambda$, over all of these
dictionaries. The fit of this function to the experimental data for
the dictionary of 64 filters is shown in~\fig{cbpdnexp09xi}, and a
corresponding performance comparison in terms of mean iterations
averaged over $\rho^{(0)}$ is displayed in~\fig{cbpdnexp09it}. Note
that none of the fixed choices of $\xi$ provide good performance over
the entire range of $\lambda$ values, while $\xi$ chosen according to
$f(\lambda)$ gives the same performance as the best choices of $\xi$
at each $\lambda$.

Additional experiments using different test images (the ``Kiel'' and
``Bridge'' standard images) as well as different dictionary filters
sizes ($12 \times 12$) indicate that $f(\lambda)$ provides a good
choice of $\xi$ over a wide range of conditions. While the choice of
$a$ giving the best fit does vary with test images, filter size, and
number of filters\footnote{The importance of selecting $\xi > 1$ for
  larger $\lambda$ values appears to be related to dictionary
  overcompleteness, corresponding to the number of filters for the
  CBPDN problem. It is also the case for the standard BPDN problem
  that the best choice of $\xi$ is greater than unity for larger
  $\lambda$ values, but for the much lower overcompleteness ratios
  usually encountered in this problem variant, the performance effect
  is far smaller, and the loss in choosing fixed $\xi = 1$ is
  usually negligible.}, the performance is not highly sensitive to
the choice of $a$ (note that the mean iteration surface for large
$\lambda$ is flat over a wide range of $\xi$ values
in~\fig{cbpdnexp02itlmdxi}) and the choice of $a = 18.3$ used
in~\fig{cbpdnexp09xi} was found to give performance at or close to the
best choice of $\xi$ in all the cases considered.

\section{Conclusion}
\label{sec:conc}

The scaling properties of the standard definitions of the primal and
dual residuals are shown to represent a potentially serious weakness
in a popular adaptive penalty strategy~\cite{he-2000-alternating} for
ADMM algorithms. The proposed solution is to normalise these residuals
so that they become invariant to scalings of the ADMM problem to which
the solution is also invariant. The impact of this issue is
demonstrated using BPDN sparse coding as an example problem. These
experiments show that the standard adaptive penalty
strategy~\cite{he-2000-alternating} performs very poorly in certain
cases, while the proposed modification based on normalised residuals is
more robust.

There is, however, a more serious issue that is not so easily
resolved: the unknown scaling relationship between the residuals and
the solution distance from optimality implies that the correct
residual ratio to target, $\xi$, is unknown, and not necessarily
unity. In some cases it is possible to construct a heuristic estimate
of this value, but it is yet to be demonstrated that such an approach
offers any real benefit over directly estimating a suitable choice of
a fixed $\rho$ parameter.

In the interests of reproducible research, software implementations of
the main algorithms proposed here are made publicly
available~\cite{wohlberg-2016-sporco}.

\appendices

\section{Scaling of the Graph Form Problem}

Many signal and image processing inverse problems can be expressed in
terms of the \emph{graph form} problem~\cite{parikh-2014-block}
\begin{equation}
\argmin_{\mb{x},\mb{z}} f(\mb{x}) + g(\mb{z}) \; \text{ such that } \;
A \mb{x} = \mb{z} \;,
\label{eq:admmgraphprob}
\end{equation}
which is a special case of~\eq{admmprob} with $B = -I$ and $\mb{c} =
0$. In this case there is a slightly different set of scalings of the
problem under which the solution is invariant, for which the most
general scaled problem $\tilde{P}$ is
\begin{equation}
  \argmin_{\mb{x},{\mb{z}}} \alpha f(\gamma
  \mb{x}) + \alpha g(\delta \mb{z}) \; \text{ s.t. } \;
  A \gamma \mb{x} = \delta \mb{z}  \;,
\label{eq:admmgrphp1}
\end{equation}
which can be expressed as graph form
problem in standard form as
\begin{equation}
  \argmin_{\mb{x},{\mb{z}}} \tl{f}(
  \mb{x}) + \tl{g}(\mb{z}) \; \text{ such that } \;
  \tl{A} \mb{x} = \mb{z}
\label{eq:admmgrphp1std}
\end{equation}
with
\begin{gather}
  \tl{f}(\mb{x}) = \alpha f(\gamma \mb{x}) \quad
  \tl{g}(\mb{z}) = \alpha g(\delta \mb{z}) \quad
  \tl{A} = \delta^{-1} \gamma A \;.
\end{gather}

The Lagrangian is
\begin{align}
  \tilde{L}({\mb{x}}, {\mb{z}}, {\mb{y}}) =
  \alpha f(\gamma \mb{x}) & + \alpha g(\delta
  \mb{z})  + \mb{y}^T (\delta^{-1} \gamma A \mb{x} - \mb{z})
\; ,
\end{align}
and the primal and dual feasibility conditions are
\begin{equation}
  \delta^{-1} \gamma A \tmb{x}^* - \tmb{z}^*  = 0 \;,
\end{equation}
and
\smallmath{
\begin{align}
  0 \in \partial \tilde{L}(\cdot, \tilde{\mb{z}}^*, \tilde{\mb{y}}^*)
  \; \Rightarrow \; & 0\!\in\! \alpha \gamma [\partial f(\cdot)](\gamma
  \tilde{\mb{x}}^*)
  +  \delta^{-1} \gamma A^T \tilde{\mb{y}}^* = 0   \\
  0 \in \partial \tilde{L}(\tilde{\mb{x}}^*, \cdot, \tilde{\mb{y}}^*)
 \; \Rightarrow \; & 0\! \in\! \alpha \delta [\partial g(\cdot)](\delta
  \tilde{\mb{z}}^*) - \tilde{\mb{y}}^* = 0
\end{align}
}
respectively.
It is easily verified that if $\mb{x}^*$, $\mb{z}^*$, and $\mb{y}^*$
satisfy the optimality criteria \eq{admmlgprim}, \eqc{admmlgrduafsg},
and \eqc{admmlgrduagsg} for problem $P$, then
\begin{gather}
\tilde{\mb{x}}^* = \gamma^{-1} \mb{x}^* \quad \;\;
\tilde{\mb{z}}^* = \delta^{-1} \mb{z}^* \quad \;\;
\tilde{\mb{y}}^* = \alpha \delta \mb{y}^*
\label{eq:optscalinggrph}
\end{gather}
satisfy the primal and dual feasibility criteria for $\tilde{P}$.

The augmented Lagrangian for $\tilde{P}$ is
\smallmath{
\begin{align}
  \tilde{L}_{\tilde{\rho}}(\mb{x}, \mb{z}, \mb{y})
  & =
\alpha f(\gamma
  \mb{x}) + \alpha g(\delta \mb{z}) \nonumber \\ &+ \alpha
  \left(\alpha^{-1} \delta^{-1} \mb{y}^T \right) ( \gamma A \mb{x} -
  \delta \mb{z} ) \nonumber \\ &+ \alpha \left(\alpha^{-1} \delta^{-2}
    \tilde{\rho} \right) \frac{1}{2} \norm{ \gamma A \mb{x} -
    \delta \mb{z}}_2^2 \;,
\end{align}
}
so that setting
\begin{equation}
 \tilde{\rho} = \alpha \delta^2 \rho
\label{eq:rhoscalegrph}
\end{equation}
gives
\begin{equation}
  \tilde{L}_{\tilde{\rho}}(\mb{x}, \mb{z}, \mb{y}) = \alpha
  L_{\rho}\left(\gamma \mb{x}, \delta \mb{z},  \alpha^{-1} \delta^{-1}
    \mb{y}\right)
  \;.
\label{eq:tlscalegrph}
\end{equation}

\section{BPDN Scaling Properties}

The scaling properties of the BPDN problem with respect to the scalar
multiplication of the input signal $\bvsigma$ depend on whether the
dictionary is considered to have fixed scaling or scale with the
signal. The former is the more common situation since the dictionary
is usually normalised, but the latter situation does occur in an
\emph{endogenous} sparse representation~\cite{dyer-2013-greedy}, in
which the signal is also used as the dictionary (with constraints on
the sparse representation to avoid the trivial solution), usually
without normalisation of the dictionary.

\subsection{Fixed Dictionary}

First, define problem $\tl{P}$ with signal $\bvsigma$ scaled by $\delta$
\begin{equation}
\argmin_{\mb{x}} \frac{1}{2} \norm{D \mb{x} - \delta \bvsigma}_2^2 +
\delta \lambda
\norm{\mb{z}}_1 \text{ s.t. } \mb{x} = \mb{z}  \;,
\end{equation}
representing the most common case in which the columns of $D$ are
normalised and $D$ does not scale with $\bvsigma$.
The corresponding Lagrangian is
\begin{align}
  \tl{L}(\mb{x}, \mb{z}, \mb{y}) &= \frac{1}{2} \norm{D \mb{x} -
    \delta \bvsigma}_2^2 + \delta \lambda
  \norm{\mb{z}}_1 + \mb{y}^T (\mb{x} - \mb{z}) \nonumber \\
  &=  \frac{1}{2} \norm{D \delta \delta^{-1} \mb{x} -
    \delta \bvsigma}_2^2 + \delta \lambda
  \norm{\delta \delta^{-1} \mb{z}}_1  \nonumber \\ & \hspace{9.5em} +
  \mb{y}^T (\delta
  \delta^{-1}\mb{x} - \delta \delta^{-1}\mb{z}) \nonumber \\
 &= \delta^2 L(\delta^{-1}\mb{x}, \delta^{-1}\mb{z},
 \delta^{-1}\mb{y}) \;.
\end{align}
Comparing with~\eq{tlscale} it is clear that we need to set
\begin{gather}
  \alpha = \delta^2 \quad \gamma = \delta^{-1} \quad \beta = \delta
\label{eq:bpdnadmmscale1}
\end{gather}
to use the ADMM scaling results of~\sctn{admmscale}. In this case the
scaling behaviour is such that changing $\delta$ does \emph{not} alter
the ratio of primal and dual residuals. Note that this merely implies
that the adaptive penalty parameter policy with standard residuals is
not \emph{guaranteed} to fail when the signal is scaled; it does not
follow that the problem scaling is such that normalised residuals are not
necessary.

\subsection{Dictionary Scales with Signal}

In the second form of scaling, $D$ is not normalised, and scales
linearly with $\bvsigma$. In this case problem $\tl{P}$ with signal
$\bvsigma$ and dictionary $D$ scaled by $\delta$ is
\begin{equation}
  \argmin_{\mb{x}} \frac{1}{2} \norm{\delta D \mb{x} - \delta \bvsigma}_2^2 +
  \delta^2 \lambda \norm{\mb{z}}_1 \text{ s.t. } \mb{x} = \mb{z}  \;.
\end{equation}
The corresponding Lagrangian is
\begin{align}
  \tl{L}(\mb{x}, \mb{z}, \mb{y}) &= \frac{1}{2} \norm{\delta D \mb{x} -
    \delta \bvsigma}_2^2 + \delta^2 \lambda
  \norm{\mb{z}}_1 + \mb{y}^T (\mb{x} - \mb{z}) \nonumber \\
  &=  \frac{\delta^2}{2} \norm{D \mb{x} -
    \bvsigma}_2^2 + \delta^2 \lambda
  \norm{\mb{z}}_1  \nonumber \\ & \hspace{7.3em} + \delta^2 \delta^{-2} \mb{y}^T (
  \mb{x} - \mb{z}) \nonumber \\
 &= \delta^2 L(\mb{x}, \mb{z}, \delta^{-2}\mb{y}) \;.
\end{align}
Comparing with~\eq{tlscale} it is clear that we need to set
\begin{gather}
  \alpha = \delta^2 \quad \gamma = 1 \quad \beta = 1
\end{gather}
to use the ADMM scaling results of~\sctn{admmscale}. In this case the
scaling behaviour is such that changing $\delta$ \emph{does} alter
the ratio of primal and dual residuals, and the adaptive
penalty parameter policy with standard residuals is guaranteed to
perform poorly for all but a restricted range of signal scaling values
$\delta$.

\section{A Degenerate Case}

An unusual degenerate case involving the TV-$\ell_1$
problem~\cite{alliney-1992-digital} illustrates that even the proposed
normalised definitions of residuals cannot always be applied without
analysis of the specific problem. This problem can be written as
\begin{equation}
  \argmin_{\mb{x}}  \| \mb{x} - \mb{s} \|_1 +  \lambda \| \sqrt{(G_0
    \mb{x})^2 + (G_1 \mb{x})^2}\|_1 \;,
\end{equation}
which can be expressed in standard ADMM form~\eq{admmprob}
(see~\cite[Sec. 2.4.4]{esser-2010-primal}) with
\begin{align}
\!\!\! f(\mb{x}) = 0 \;\;\;\;\; g(\mb{z}) = \| \mb{z}_s - \mb{s} \|_1 +  \lambda
\| \sqrt{\mb{z}_0^2 + \mb{z}_1^2}\|_1 \nonumber \\
A = \left( \begin{array}{c} G_0 \\ G_1 \\ I \end{array} \right) \;
B = -I \;\;\;
\mb{c} = \left( \begin{array}{c} 0 \\ 0 \\ \mb{s} \end{array} \right)
\;
\mb{z} = \left( \begin{array}{c} \mb{z}_0 \\ \mb{z}_1 \\
    \mb{z}_s \end{array} \right) .
\end{align}

Since $f(\mb{x}) = 0$, dual feasibility condition~\eq{admmlgrduafsg}
is simply $A^T \mb{y}^* = 0$ and~\eq{dualres0}, from which the
definition~\eq{dualres} is derived, degenerates to
\begin{equation}
  \rho A^T B (\mb{z}^{(k+1)} - \mb{z}^{(k)}) = A^T \mb{y}^{(k+1)} \;.
\label{eq:dualresdegen0}
\end{equation}
Clearly $A^T \mb{y}^{(k+1)}$ is unsuitable either as a normalisation
term for the dual residual or as a factor in the stopping tolerance.

\bibliographystyle{IEEEtranD}
\bibliography{admm}

\end{document}

%% file: otherexprmnt01fnval.pstex_t
% GNUPLOT: LaTeX picture with Postscript
\begingroup%
\makeatletter%
\newcommand{\GNUPLOTspecial}{%
  \@sanitize\catcode`\%=14\relax\special}%
\setlength{\unitlength}{0.0500bp}%
\begin{picture}(5616,3124)(0,0)%
  \put(4252,2521){\makebox(0,0)[r]{\strut{}Normalised}}%
  \put(4252,2721){\makebox(0,0)[r]{\strut{}Standard}}%
  \put(3127,140){\makebox(0,0){\strut{}Iteration number}}%
  \put(160,1762){%
\rotatebox{-270}{%
  \makebox(0,0){\strut{}$(p^{(k)} - p^*) / p^*$\vspace{-7mm}}%
}}%
  \put(5155,440){\makebox(0,0){\strut{} 1000}}%
  \put(4750,440){\makebox(0,0){\strut{} 900}}%
  \put(4344,440){\makebox(0,0){\strut{} 800}}%
  \put(3939,440){\makebox(0,0){\strut{} 700}}%
  \put(3533,440){\makebox(0,0){\strut{} 600}}%
  \put(3128,440){\makebox(0,0){\strut{} 500}}%
  \put(2722,440){\makebox(0,0){\strut{} 400}}%
  \put(2317,440){\makebox(0,0){\strut{} 300}}%
  \put(1911,440){\makebox(0,0){\strut{} 200}}%
  \put(1506,440){\makebox(0,0){\strut{} 100}}%
  \put(1100,440){\makebox(0,0){\strut{} 0}}%
  \put(980,2884){\makebox(0,0)[r]{\strut{}1e+02}}%
  \put(980,2635){\makebox(0,0)[r]{\strut{}1e+00}}%
  \put(980,2385){\makebox(0,0)[r]{\strut{}1e-02}}%
  \put(980,2136){\makebox(0,0)[r]{\strut{}1e-04}}%
  \put(980,1887){\makebox(0,0)[r]{\strut{}1e-06}}%
  \put(980,1637){\makebox(0,0)[r]{\strut{}1e-08}}%
  \put(980,1388){\makebox(0,0)[r]{\strut{}1e-10}}%
  \put(980,1139){\makebox(0,0)[r]{\strut{}1e-12}}%
  \put(980,889){\makebox(0,0)[r]{\strut{}1e-14}}%
  \put(980,640){\makebox(0,0)[r]{\strut{}1e-16}}%
\includegraphics{otherexprmnt01fnval}%
\end{picture}%
\endgroup
 

%% file: otherexprmnt01priduares.pstex_t
% GNUPLOT: LaTeX picture with Postscript
\begingroup%
\makeatletter%
\newcommand{\GNUPLOTspecial}{%
  \@sanitize\catcode`\%=14\relax\special}%
\setlength{\unitlength}{0.0500bp}%
\begin{picture}(5616,3124)(0,0)%
  \put(4612,1311){\makebox(0,0)[r]{\strut{}Dual Normalised}}%
  \put(4612,1511){\makebox(0,0)[r]{\strut{}Primal Normalised}}%
  \put(4612,1711){\makebox(0,0)[r]{\strut{}Dual Standard}}%
  \put(4612,1911){\makebox(0,0)[r]{\strut{}Primal Standard}}%
  \put(3247,140){\makebox(0,0){\strut{}Iteration number}}%
  \put(160,1762){%
\rotatebox{-270}{%
  \makebox(0,0){\strut{}Residual\vspace{-7mm}}%
}}%
  \put(5155,440){\makebox(0,0){\strut{} 1000}}%
  \put(4774,440){\makebox(0,0){\strut{} 900}}%
  \put(4392,440){\makebox(0,0){\strut{} 800}}%
  \put(4011,440){\makebox(0,0){\strut{} 700}}%
  \put(3629,440){\makebox(0,0){\strut{} 600}}%
  \put(3248,440){\makebox(0,0){\strut{} 500}}%
  \put(2866,440){\makebox(0,0){\strut{} 400}}%
  \put(2485,440){\makebox(0,0){\strut{} 300}}%
  \put(2103,440){\makebox(0,0){\strut{} 200}}%
  \put(1722,440){\makebox(0,0){\strut{} 100}}%
  \put(1340,440){\makebox(0,0){\strut{} 0}}%
  \put(1220,2884){\makebox(0,0)[r]{\strut{}1.0e+04}}%
  \put(1220,2635){\makebox(0,0)[r]{\strut{}1.0e+02}}%
  \put(1220,2385){\makebox(0,0)[r]{\strut{}1.0e+00}}%
  \put(1220,2136){\makebox(0,0)[r]{\strut{}1.0e-02}}%
  \put(1220,1887){\makebox(0,0)[r]{\strut{}1.0e-04}}%
  \put(1220,1637){\makebox(0,0)[r]{\strut{}1.0e-06}}%
  \put(1220,1388){\makebox(0,0)[r]{\strut{}1.0e-08}}%
  \put(1220,1139){\makebox(0,0)[r]{\strut{}1.0e-10}}%
  \put(1220,889){\makebox(0,0)[r]{\strut{}1.0e-12}}%
  \put(1220,640){\makebox(0,0)[r]{\strut{}1.0e-14}}%
\includegraphics{otherexprmnt01priduares}%
\end{picture}%
\endgroup
 

%% file: otherexprmnt01rho.pstex_t
% GNUPLOT: LaTeX picture with Postscript
\begingroup%
\makeatletter%
\newcommand{\GNUPLOTspecial}{%
  \@sanitize\catcode`\%=14\relax\special}%
\setlength{\unitlength}{0.0500bp}%
\begin{picture}(5616,3124)(0,0)%
  \put(4252,2521){\makebox(0,0)[r]{\strut{}Normalised}}%
  \put(4252,2721){\makebox(0,0)[r]{\strut{}Standard}}%
  \put(3247,140){\makebox(0,0){\strut{}Iteration number}}%
  \put(160,1762){%
\rotatebox{-270}{%
  \makebox(0,0){\strut{}$\rho$\vspace{-7mm}}%
}}%
  \put(5155,440){\makebox(0,0){\strut{} 1000}}%
  \put(4774,440){\makebox(0,0){\strut{} 900}}%
  \put(4392,440){\makebox(0,0){\strut{} 800}}%
  \put(4011,440){\makebox(0,0){\strut{} 700}}%
  \put(3629,440){\makebox(0,0){\strut{} 600}}%
  \put(3248,440){\makebox(0,0){\strut{} 500}}%
  \put(2866,440){\makebox(0,0){\strut{} 400}}%
  \put(2485,440){\makebox(0,0){\strut{} 300}}%
  \put(2103,440){\makebox(0,0){\strut{} 200}}%
  \put(1722,440){\makebox(0,0){\strut{} 100}}%
  \put(1340,440){\makebox(0,0){\strut{} 0}}%
  \put(1220,2884){\makebox(0,0)[r]{\strut{}1.2e+02}}%
  \put(1220,2510){\makebox(0,0)[r]{\strut{}1.0e+02}}%
  \put(1220,2136){\makebox(0,0)[r]{\strut{}8.0e+01}}%
  \put(1220,1762){\makebox(0,0)[r]{\strut{}6.0e+01}}%
  \put(1220,1388){\makebox(0,0)[r]{\strut{}4.0e+01}}%
  \put(1220,1014){\makebox(0,0)[r]{\strut{}2.0e+01}}%
  \put(1220,640){\makebox(0,0)[r]{\strut{}0.0e+00}}%
\includegraphics{otherexprmnt01rho}%
\end{picture}%
\endgroup
 

%% file: expbpdn12_itrh_std_1en03_128.pstex_t
% GNUPLOT: LaTeX picture with Postscript
\begingroup%
\makeatletter%
\newcommand{\GNUPLOTspecial}{%
  \@sanitize\catcode`\%=14\relax\special}%
\setlength{\unitlength}{0.0500bp}%
\begin{picture}(5615,3124)(0,0)%
  \put(4492,1721){\makebox(0,0)[r]{\strut{}1.2/Auto}}%
  \put(4492,1921){\makebox(0,0)[r]{\strut{}2/Auto}}%
  \put(4492,2121){\makebox(0,0)[r]{\strut{}2/1.2}}%
  \put(4492,2321){\makebox(0,0)[r]{\strut{}5/2}}%
  \put(4492,2521){\makebox(0,0)[r]{\strut{}10/2}}%
  \put(4492,2721){\makebox(0,0)[r]{\strut{}Fixed}}%
  \put(2967,140){\makebox(0,0){\strut{}$\rho^{(0)}$}}%
  \put(224,1762){%
\rotatebox{-270}{%
  \makebox(0,0){\strut{}Iterations}%
}}%
  \put(5155,440){\makebox(0,0){\strut{}1e+00}}%
  \put(4061,440){\makebox(0,0){\strut{}1e-01}}%
  \put(2967,440){\makebox(0,0){\strut{}1e-02}}%
  \put(1874,440){\makebox(0,0){\strut{}1e-03}}%
  \put(780,440){\makebox(0,0){\strut{}1e-04}}%
  \put(660,2884){\makebox(0,0)[r]{\strut{} 500}}%
  \put(660,2435){\makebox(0,0)[r]{\strut{} 400}}%
  \put(660,1986){\makebox(0,0)[r]{\strut{} 300}}%
  \put(660,1538){\makebox(0,0)[r]{\strut{} 200}}%
  \put(660,1089){\makebox(0,0)[r]{\strut{} 100}}%
  \put(660,640){\makebox(0,0)[r]{\strut{} 0}}%
\includegraphics{expbpdn12_itrh_std_1en03_128}%
\end{picture}%
\endgroup
 

%% file: expbpdn12_itrh_nrm_1en03_128.pstex_t
% GNUPLOT: LaTeX picture with Postscript
\begingroup%
\makeatletter%
\newcommand{\GNUPLOTspecial}{%
  \@sanitize\catcode`\%=14\relax\special}%
\setlength{\unitlength}{0.0500bp}%
\begin{picture}(5615,3124)(0,0)%
  \put(3056,1721){\makebox(0,0)[r]{\strut{}1.2/Auto}}%
  \put(3056,1921){\makebox(0,0)[r]{\strut{}2/Auto}}%
  \put(3056,2121){\makebox(0,0)[r]{\strut{}2/1.2}}%
  \put(3056,2321){\makebox(0,0)[r]{\strut{}5/2}}%
  \put(3056,2521){\makebox(0,0)[r]{\strut{}10/2}}%
  \put(3056,2721){\makebox(0,0)[r]{\strut{}Fixed}}%
  \put(2967,140){\makebox(0,0){\strut{}$\rho^{(0)}$}}%
  \put(224,1762){%
\rotatebox{-270}{%
  \makebox(0,0){\strut{}Iterations}%
}}%
  \put(5155,440){\makebox(0,0){\strut{}1e+00}}%
  \put(4061,440){\makebox(0,0){\strut{}1e-01}}%
  \put(2967,440){\makebox(0,0){\strut{}1e-02}}%
  \put(1874,440){\makebox(0,0){\strut{}1e-03}}%
  \put(780,440){\makebox(0,0){\strut{}1e-04}}%
  \put(660,2884){\makebox(0,0)[r]{\strut{} 500}}%
  \put(660,2435){\makebox(0,0)[r]{\strut{} 400}}%
  \put(660,1986){\makebox(0,0)[r]{\strut{} 300}}%
  \put(660,1538){\makebox(0,0)[r]{\strut{} 200}}%
  \put(660,1089){\makebox(0,0)[r]{\strut{} 100}}%
  \put(660,640){\makebox(0,0)[r]{\strut{} 0}}%
\includegraphics{expbpdn12_itrh_nrm_1en03_128}%
\end{picture}%
\endgroup
 

%% file: expbpdn12_itlm_std_64.pstex_t
% GNUPLOT: LaTeX picture with Postscript
\begingroup%
\makeatletter%
\newcommand{\GNUPLOTspecial}{%
  \@sanitize\catcode`\%=14\relax\special}%
\setlength{\unitlength}{0.0500bp}%
\begin{picture}(5616,3528)(0,0)%
  \put(4710,3275){\makebox(0,0)[r]{\strut{}1.2/Auto}}%
  \put(3567,3075){\makebox(0,0)[r]{\strut{}2/Auto}}%
  \put(3567,3275){\makebox(0,0)[r]{\strut{}2/1.2}}%
  \put(2424,3075){\makebox(0,0)[r]{\strut{}5/2}}%
  \put(2424,3275){\makebox(0,0)[r]{\strut{}10/2}}%
  \put(1281,3075){\makebox(0,0)[r]{\strut{}Fixed}}%
  \put(1281,3275){\makebox(0,0)[r]{\strut{}Fixed (min)}}%
  \put(2967,140){\makebox(0,0){\strut{}$\lambda$}}%
  \put(224,1753){%
\rotatebox{-270}{%
  \makebox(0,0){\strut{}Mean Iterations}%
}}%
  \put(5155,440){\makebox(0,0){\strut{}3e-01}}%
  \put(4459,440){\makebox(0,0){\strut{}1e-01}}%
  \put(3697,440){\makebox(0,0){\strut{}3e-02}}%
  \put(3001,440){\makebox(0,0){\strut{}1e-02}}%
  \put(2238,440){\makebox(0,0){\strut{}3e-03}}%
  \put(1543,440){\makebox(0,0){\strut{}1e-03}}%
  \put(780,440){\makebox(0,0){\strut{}3e-04}}%
  \put(660,2867){\makebox(0,0)[r]{\strut{} 500}}%
  \put(660,2422){\makebox(0,0)[r]{\strut{} 400}}%
  \put(660,1976){\makebox(0,0)[r]{\strut{} 300}}%
  \put(660,1531){\makebox(0,0)[r]{\strut{} 200}}%
  \put(660,1085){\makebox(0,0)[r]{\strut{} 100}}%
  \put(660,640){\makebox(0,0)[r]{\strut{} 0}}%
\includegraphics{expbpdn12_itlm_std_64}%
\end{picture}%
\endgroup
 

%% file: expbpdn12_itlm_nrm_64.pstex_t
% GNUPLOT: LaTeX picture with Postscript
\begingroup%
\makeatletter%
\newcommand{\GNUPLOTspecial}{%
  \@sanitize\catcode`\%=14\relax\special}%
\setlength{\unitlength}{0.0500bp}%
\begin{picture}(5616,3528)(0,0)%
  \put(4710,3275){\makebox(0,0)[r]{\strut{}1.2/Auto}}%
  \put(3567,3075){\makebox(0,0)[r]{\strut{}2/Auto}}%
  \put(3567,3275){\makebox(0,0)[r]{\strut{}2/1.2}}%
  \put(2424,3075){\makebox(0,0)[r]{\strut{}5/2}}%
  \put(2424,3275){\makebox(0,0)[r]{\strut{}10/2}}%
  \put(1281,3075){\makebox(0,0)[r]{\strut{}Fixed}}%
  \put(1281,3275){\makebox(0,0)[r]{\strut{}Fixed (min)}}%
  \put(2967,140){\makebox(0,0){\strut{}$\lambda$}}%
  \put(224,1753){%
\rotatebox{-270}{%
  \makebox(0,0){\strut{}Mean Iterations}%
}}%
  \put(5155,440){\makebox(0,0){\strut{}3e-01}}%
  \put(4459,440){\makebox(0,0){\strut{}1e-01}}%
  \put(3697,440){\makebox(0,0){\strut{}3e-02}}%
  \put(3001,440){\makebox(0,0){\strut{}1e-02}}%
  \put(2238,440){\makebox(0,0){\strut{}3e-03}}%
  \put(1543,440){\makebox(0,0){\strut{}1e-03}}%
  \put(780,440){\makebox(0,0){\strut{}3e-04}}%
  \put(660,2867){\makebox(0,0)[r]{\strut{} 500}}%
  \put(660,2422){\makebox(0,0)[r]{\strut{} 400}}%
  \put(660,1976){\makebox(0,0)[r]{\strut{} 300}}%
  \put(660,1531){\makebox(0,0)[r]{\strut{} 200}}%
  \put(660,1085){\makebox(0,0)[r]{\strut{} 100}}%
  \put(660,640){\makebox(0,0)[r]{\strut{} 0}}%
\includegraphics{expbpdn12_itlm_nrm_64}%
\end{picture}%
\endgroup
 

%% file: expbpdn12_itlm_std_96.pstex_t
% GNUPLOT: LaTeX picture with Postscript
\begingroup%
\makeatletter%
\newcommand{\GNUPLOTspecial}{%
  \@sanitize\catcode`\%=14\relax\special}%
\setlength{\unitlength}{0.0500bp}%
\begin{picture}(5615,3527)(0,0)%
  \put(4710,3274){\makebox(0,0)[r]{\strut{}1.2/Auto}}%
  \put(3567,3074){\makebox(0,0)[r]{\strut{}2/Auto}}%
  \put(3567,3274){\makebox(0,0)[r]{\strut{}2/1.2}}%
  \put(2424,3074){\makebox(0,0)[r]{\strut{}5/2}}%
  \put(2424,3274){\makebox(0,0)[r]{\strut{}10/2}}%
  \put(1281,3074){\makebox(0,0)[r]{\strut{}Fixed}}%
  \put(1281,3274){\makebox(0,0)[r]{\strut{}Fixed (min)}}%
  \put(2967,140){\makebox(0,0){\strut{}$\lambda$}}%
  \put(224,1753){%
\rotatebox{-270}{%
  \makebox(0,0){\strut{}Mean Iterations}%
}}%
  \put(5155,440){\makebox(0,0){\strut{}3e-01}}%
  \put(4459,440){\makebox(0,0){\strut{}1e-01}}%
  \put(3697,440){\makebox(0,0){\strut{}3e-02}}%
  \put(3001,440){\makebox(0,0){\strut{}1e-02}}%
  \put(2238,440){\makebox(0,0){\strut{}3e-03}}%
  \put(1543,440){\makebox(0,0){\strut{}1e-03}}%
  \put(780,440){\makebox(0,0){\strut{}3e-04}}%
  \put(660,2867){\makebox(0,0)[r]{\strut{} 500}}%
  \put(660,2422){\makebox(0,0)[r]{\strut{} 400}}%
  \put(660,1976){\makebox(0,0)[r]{\strut{} 300}}%
  \put(660,1531){\makebox(0,0)[r]{\strut{} 200}}%
  \put(660,1085){\makebox(0,0)[r]{\strut{} 100}}%
  \put(660,640){\makebox(0,0)[r]{\strut{} 0}}%
\includegraphics{expbpdn12_itlm_std_96}%
\end{picture}%
\endgroup
 

%% file: expbpdn12_itlm_nrm_96.pstex_t
% GNUPLOT: LaTeX picture with Postscript
\begingroup%
\makeatletter%
\newcommand{\GNUPLOTspecial}{%
  \@sanitize\catcode`\%=14\relax\special}%
\setlength{\unitlength}{0.0500bp}%
\begin{picture}(5615,3527)(0,0)%
  \put(4710,3274){\makebox(0,0)[r]{\strut{}1.2/Auto}}%
  \put(3567,3074){\makebox(0,0)[r]{\strut{}2/Auto}}%
  \put(3567,3274){\makebox(0,0)[r]{\strut{}2/1.2}}%
  \put(2424,3074){\makebox(0,0)[r]{\strut{}5/2}}%
  \put(2424,3274){\makebox(0,0)[r]{\strut{}10/2}}%
  \put(1281,3074){\makebox(0,0)[r]{\strut{}Fixed}}%
  \put(1281,3274){\makebox(0,0)[r]{\strut{}Fixed (min)}}%
  \put(2967,140){\makebox(0,0){\strut{}$\lambda$}}%
  \put(224,1753){%
\rotatebox{-270}{%
  \makebox(0,0){\strut{}Mean Iterations}%
}}%
  \put(5155,440){\makebox(0,0){\strut{}3e-01}}%
  \put(4459,440){\makebox(0,0){\strut{}1e-01}}%
  \put(3697,440){\makebox(0,0){\strut{}3e-02}}%
  \put(3001,440){\makebox(0,0){\strut{}1e-02}}%
  \put(2238,440){\makebox(0,0){\strut{}3e-03}}%
  \put(1543,440){\makebox(0,0){\strut{}1e-03}}%
  \put(780,440){\makebox(0,0){\strut{}3e-04}}%
  \put(660,2867){\makebox(0,0)[r]{\strut{} 500}}%
  \put(660,2422){\makebox(0,0)[r]{\strut{} 400}}%
  \put(660,1976){\makebox(0,0)[r]{\strut{} 300}}%
  \put(660,1531){\makebox(0,0)[r]{\strut{} 200}}%
  \put(660,1085){\makebox(0,0)[r]{\strut{} 100}}%
  \put(660,640){\makebox(0,0)[r]{\strut{} 0}}%
\includegraphics{expbpdn12_itlm_nrm_96}%
\end{picture}%
\endgroup
 

%% file: expbpdn12_itlm_std_128.pstex_t
% GNUPLOT: LaTeX picture with Postscript
\begingroup%
\makeatletter%
\newcommand{\GNUPLOTspecial}{%
  \@sanitize\catcode`\%=14\relax\special}%
\setlength{\unitlength}{0.0500bp}%
\begin{picture}(5615,3527)(0,0)%
  \put(4710,3274){\makebox(0,0)[r]{\strut{}1.2/Auto}}%
  \put(3567,3074){\makebox(0,0)[r]{\strut{}2/Auto}}%
  \put(3567,3274){\makebox(0,0)[r]{\strut{}2/1.2}}%
  \put(2424,3074){\makebox(0,0)[r]{\strut{}5/2}}%
  \put(2424,3274){\makebox(0,0)[r]{\strut{}10/2}}%
  \put(1281,3074){\makebox(0,0)[r]{\strut{}Fixed}}%
  \put(1281,3274){\makebox(0,0)[r]{\strut{}Fixed (min)}}%
  \put(2967,140){\makebox(0,0){\strut{}$\lambda$}}%
  \put(224,1753){%
\rotatebox{-270}{%
  \makebox(0,0){\strut{}Mean Iterations}%
}}%
  \put(5155,440){\makebox(0,0){\strut{}3e-01}}%
  \put(4459,440){\makebox(0,0){\strut{}1e-01}}%
  \put(3697,440){\makebox(0,0){\strut{}3e-02}}%
  \put(3001,440){\makebox(0,0){\strut{}1e-02}}%
  \put(2238,440){\makebox(0,0){\strut{}3e-03}}%
  \put(1543,440){\makebox(0,0){\strut{}1e-03}}%
  \put(780,440){\makebox(0,0){\strut{}3e-04}}%
  \put(660,2867){\makebox(0,0)[r]{\strut{} 500}}%
  \put(660,2422){\makebox(0,0)[r]{\strut{} 400}}%
  \put(660,1976){\makebox(0,0)[r]{\strut{} 300}}%
  \put(660,1531){\makebox(0,0)[r]{\strut{} 200}}%
  \put(660,1085){\makebox(0,0)[r]{\strut{} 100}}%
  \put(660,640){\makebox(0,0)[r]{\strut{} 0}}%
\includegraphics{expbpdn12_itlm_std_128}%
\end{picture}%
\endgroup
 

%% file: expbpdn12_itlm_nrm_128.pstex_t
% GNUPLOT: LaTeX picture with Postscript
\begingroup%
\makeatletter%
\newcommand{\GNUPLOTspecial}{%
  \@sanitize\catcode`\%=14\relax\special}%
\setlength{\unitlength}{0.0500bp}%
\begin{picture}(5615,3527)(0,0)%
  \put(4710,3274){\makebox(0,0)[r]{\strut{}1.2/Auto}}%
  \put(3567,3074){\makebox(0,0)[r]{\strut{}2/Auto}}%
  \put(3567,3274){\makebox(0,0)[r]{\strut{}2/1.2}}%
  \put(2424,3074){\makebox(0,0)[r]{\strut{}5/2}}%
  \put(2424,3274){\makebox(0,0)[r]{\strut{}10/2}}%
  \put(1281,3074){\makebox(0,0)[r]{\strut{}Fixed}}%
  \put(1281,3274){\makebox(0,0)[r]{\strut{}Fixed (min)}}%
  \put(2967,140){\makebox(0,0){\strut{}$\lambda$}}%
  \put(224,1753){%
\rotatebox{-270}{%
  \makebox(0,0){\strut{}Mean Iterations}%
}}%
  \put(5155,440){\makebox(0,0){\strut{}3e-01}}%
  \put(4459,440){\makebox(0,0){\strut{}1e-01}}%
  \put(3697,440){\makebox(0,0){\strut{}3e-02}}%
  \put(3001,440){\makebox(0,0){\strut{}1e-02}}%
  \put(2238,440){\makebox(0,0){\strut{}3e-03}}%
  \put(1543,440){\makebox(0,0){\strut{}1e-03}}%
  \put(780,440){\makebox(0,0){\strut{}3e-04}}%
  \put(660,2867){\makebox(0,0)[r]{\strut{} 500}}%
  \put(660,2422){\makebox(0,0)[r]{\strut{} 400}}%
  \put(660,1976){\makebox(0,0)[r]{\strut{} 300}}%
  \put(660,1531){\makebox(0,0)[r]{\strut{} 200}}%
  \put(660,1085){\makebox(0,0)[r]{\strut{} 100}}%
  \put(660,640){\makebox(0,0)[r]{\strut{} 0}}%
\includegraphics{expbpdn12_itlm_nrm_128}%
\end{picture}%
\endgroup
 

%% file: expcbpdn04_itlm_std_64.pstex_t
% GNUPLOT: LaTeX picture with Postscript
\begingroup%
\makeatletter%
\newcommand{\GNUPLOTspecial}{%
  \@sanitize\catcode`\%=14\relax\special}%
\setlength{\unitlength}{0.0500bp}%
\begin{picture}(5615,3527)(0,0)%
  \put(4710,3274){\makebox(0,0)[r]{\strut{}1.2/Auto}}%
  \put(3567,3074){\makebox(0,0)[r]{\strut{}2/Auto}}%
  \put(3567,3274){\makebox(0,0)[r]{\strut{}2/1.2}}%
  \put(2424,3074){\makebox(0,0)[r]{\strut{}5/2}}%
  \put(2424,3274){\makebox(0,0)[r]{\strut{}10/2}}%
  \put(1281,3074){\makebox(0,0)[r]{\strut{}Fixed}}%
  \put(1281,3274){\makebox(0,0)[r]{\strut{}Fixed (min)}}%
  \put(2967,140){\makebox(0,0){\strut{}$\lambda$}}%
  \put(224,1753){%
\rotatebox{-270}{%
  \makebox(0,0){\strut{}Mean Iterations}%
}}%
  \put(5155,440){\makebox(0,0){\strut{}3e-01}}%
  \put(4312,440){\makebox(0,0){\strut{}1e-01}}%
  \put(3389,440){\makebox(0,0){\strut{}3e-02}}%
  \put(2546,440){\makebox(0,0){\strut{}1e-02}}%
  \put(1623,440){\makebox(0,0){\strut{}3e-03}}%
  \put(780,440){\makebox(0,0){\strut{}1e-03}}%
  \put(660,2867){\makebox(0,0)[r]{\strut{} 500}}%
  \put(660,2422){\makebox(0,0)[r]{\strut{} 400}}%
  \put(660,1976){\makebox(0,0)[r]{\strut{} 300}}%
  \put(660,1531){\makebox(0,0)[r]{\strut{} 200}}%
  \put(660,1085){\makebox(0,0)[r]{\strut{} 100}}%
  \put(660,640){\makebox(0,0)[r]{\strut{} 0}}%
\includegraphics{expcbpdn04_itlm_std_64}%
\end{picture}%
\endgroup
 

%% file: expcbpdn04_itlm_nrm_64.pstex_t
% GNUPLOT: LaTeX picture with Postscript
\begingroup%
\makeatletter%
\newcommand{\GNUPLOTspecial}{%
  \@sanitize\catcode`\%=14\relax\special}%
\setlength{\unitlength}{0.0500bp}%
\begin{picture}(5615,3527)(0,0)%
  \put(4710,3274){\makebox(0,0)[r]{\strut{}1.2/Auto}}%
  \put(3567,3074){\makebox(0,0)[r]{\strut{}2/Auto}}%
  \put(3567,3274){\makebox(0,0)[r]{\strut{}2/1.2}}%
  \put(2424,3074){\makebox(0,0)[r]{\strut{}5/2}}%
  \put(2424,3274){\makebox(0,0)[r]{\strut{}10/2}}%
  \put(1281,3074){\makebox(0,0)[r]{\strut{}Fixed}}%
  \put(1281,3274){\makebox(0,0)[r]{\strut{}Fixed (min)}}%
  \put(2967,140){\makebox(0,0){\strut{}$\lambda$}}%
  \put(224,1753){%
\rotatebox{-270}{%
  \makebox(0,0){\strut{}Mean Iterations}%
}}%
  \put(5155,440){\makebox(0,0){\strut{}3e-01}}%
  \put(4312,440){\makebox(0,0){\strut{}1e-01}}%
  \put(3389,440){\makebox(0,0){\strut{}3e-02}}%
  \put(2546,440){\makebox(0,0){\strut{}1e-02}}%
  \put(1623,440){\makebox(0,0){\strut{}3e-03}}%
  \put(780,440){\makebox(0,0){\strut{}1e-03}}%
  \put(660,2867){\makebox(0,0)[r]{\strut{} 500}}%
  \put(660,2422){\makebox(0,0)[r]{\strut{} 400}}%
  \put(660,1976){\makebox(0,0)[r]{\strut{} 300}}%
  \put(660,1531){\makebox(0,0)[r]{\strut{} 200}}%
  \put(660,1085){\makebox(0,0)[r]{\strut{} 100}}%
  \put(660,640){\makebox(0,0)[r]{\strut{} 0}}%
\includegraphics{expcbpdn04_itlm_nrm_64}%
\end{picture}%
\endgroup
 

%% file: cbpdn_exp09_rs_xmpl.pstex_t
% GNUPLOT: LaTeX picture with Postscript
\begingroup%
\makeatletter%
\newcommand{\GNUPLOTspecial}{%
  \@sanitize\catcode`\%=14\relax\special}%
\setlength{\unitlength}{0.0500bp}%
\begin{picture}(5615,3124)(0,0)%
  \put(4252,2121){\makebox(0,0)[r]{\strut{}Dual $\xi = 5.0$}}%
  \put(4252,2321){\makebox(0,0)[r]{\strut{}Primal $\xi = 5.0$}}%
  \put(4252,2521){\makebox(0,0)[r]{\strut{}Dual $\xi = 1.0$}}%
  \put(4252,2721){\makebox(0,0)[r]{\strut{}Primal $\xi = 1.0$}}%
  \put(3127,140){\makebox(0,0){\strut{}Iterations}}%
  \put(160,1762){%
\rotatebox{-270}{%
  \makebox(0,0){\strut{}Residual\vspace{-7.6mm}}%
}}%
  \put(5155,440){\makebox(0,0){\strut{} 140}}%
  \put(4576,440){\makebox(0,0){\strut{} 120}}%
  \put(3996,440){\makebox(0,0){\strut{} 100}}%
  \put(3417,440){\makebox(0,0){\strut{} 80}}%
  \put(2838,440){\makebox(0,0){\strut{} 60}}%
  \put(2259,440){\makebox(0,0){\strut{} 40}}%
  \put(1679,440){\makebox(0,0){\strut{} 20}}%
  \put(1100,440){\makebox(0,0){\strut{} 0}}%
  \put(980,2884){\makebox(0,0)[r]{\strut{}1e+00}}%
  \put(980,2323){\makebox(0,0)[r]{\strut{}1e-01}}%
  \put(980,1762){\makebox(0,0)[r]{\strut{}1e-02}}%
  \put(980,1201){\makebox(0,0)[r]{\strut{}1e-03}}%
  \put(980,640){\makebox(0,0)[r]{\strut{}1e-04}}%
\includegraphics{cbpdn_exp09_rs_xmpl}%
\end{picture}%
\endgroup
 

%% file: cbpdn_newexp02_mnitlmdxi_64.pstex_t
% GNUPLOT: LaTeX picture with Postscript
\begingroup%
\makeatletter%
\newcommand{\GNUPLOTspecial}{%
  \@sanitize\catcode`\%=14\relax\special}%
\setlength{\unitlength}{0.0500bp}%
\begin{picture}(5400,3780)(0,0)%
  \put(-62,2028){%
\rotatebox{-270}{%
  \makebox(0,0){\strut{}Mean iterations\vspace{-7mm}}%
}}%
  \put(652,2461){\makebox(0,0)[r]{\strut{} 500}}%
  \put(652,2244){\makebox(0,0)[r]{\strut{} 400}}%
  \put(652,2026){\makebox(0,0)[r]{\strut{} 300}}%
  \put(652,1809){\makebox(0,0)[r]{\strut{} 200}}%
  \put(1436,564){\makebox(0,0){\strut{}$\xi$}}%
  \put(2862,450){\makebox(0,0){\strut{} 10}}%
  \put(2436,561){\makebox(0,0){\strut{} 5}}%
  \put(2122,643){\makebox(0,0){\strut{} 3}}%
  \put(1445,818){\makebox(0,0){\strut{} 1}}%
  \put(1019,929){\makebox(0,0){\strut{} 0.5}}%
  \put(704,1011){\makebox(0,0){\strut{} 0.3}}%
  \put(4318,701){\makebox(0,0){\strut{}$\lambda$}}%
  \put(3040,531){\makebox(0,0)[l]{\strut{} 0.3}}%
  \put(3365,669){\makebox(0,0)[l]{\strut{} 0.1}}%
  \put(3720,821){\makebox(0,0)[l]{\strut{} 0.03}}%
  \put(4045,959){\makebox(0,0)[l]{\strut{} 0.01}}%
  \put(4401,1110){\makebox(0,0)[l]{\strut{} 0.003}}%
  \put(4726,1248){\makebox(0,0)[l]{\strut{} 0.001}}%
\includegraphics{cbpdn_newexp02_mnitlmdxi_64}%
\end{picture}%
\endgroup
 

%% file: cbpdn_newexp02_sditlmdxi_64.pstex_t
% GNUPLOT: LaTeX picture with Postscript
\begingroup%
\makeatletter%
\newcommand{\GNUPLOTspecial}{%
  \@sanitize\catcode`\%=14\relax\special}%
\setlength{\unitlength}{0.0500bp}%
\begin{picture}(5400,3780)(0,0)%
  \put(-41,1901){%
\rotatebox{-270}{%
  \makebox(0,0){\strut{}Std. dev. of iter.\vspace{-7mm}}%
}}%
  \put(673,2052){\makebox(0,0)[r]{\strut{} 30}}%
  \put(673,1887){\makebox(0,0)[r]{\strut{} 20}}%
  \put(673,1721){\makebox(0,0)[r]{\strut{} 10}}%
  \put(673,1555){\makebox(0,0)[r]{\strut{} 0}}%
  \put(1552,486){\makebox(0,0){\strut{}$\xi$}}%
  \put(3002,367){\makebox(0,0){\strut{} 10}}%
  \put(2554,499){\makebox(0,0){\strut{} 5}}%
  \put(2223,597){\makebox(0,0){\strut{} 3}}%
  \put(1512,807){\makebox(0,0){\strut{} 1}}%
  \put(1063,939){\makebox(0,0){\strut{} 0.5}}%
  \put(732,1037){\makebox(0,0){\strut{} 0.3}}%
  \put(4403,749){\makebox(0,0){\strut{}$\lambda$}}%
  \put(3182,454){\makebox(0,0)[l]{\strut{} 0.3}}%
  \put(3477,645){\makebox(0,0)[l]{\strut{} 0.1}}%
  \put(3800,855){\makebox(0,0)[l]{\strut{} 0.03}}%
  \put(4095,1046){\makebox(0,0)[l]{\strut{} 0.01}}%
  \put(4418,1256){\makebox(0,0)[l]{\strut{} 0.003}}%
  \put(4713,1447){\makebox(0,0)[l]{\strut{} 0.001}}%
\includegraphics{cbpdn_newexp02_sditlmdxi_64}%
\end{picture}%
\endgroup
 

%% file: cbpdn_newexp02_lmdxi_64.pstex_t
% GNUPLOT: LaTeX picture with Postscript
\begingroup%
\makeatletter%
\newcommand{\GNUPLOTspecial}{%
  \@sanitize\catcode`\%=14\relax\special}%
\setlength{\unitlength}{0.0500bp}%
\begin{picture}(5615,3124)(0,0)%
  \put(2180,2521){\makebox(0,0)[r]{\strut{}$f(\lambda)$}}%
  \put(2180,2721){\makebox(0,0)[r]{\strut{}Experiment}}%
  \put(3007,140){\makebox(0,0){\strut{}$\lambda$}}%
  \put(160,1762){%
\rotatebox{-270}{%
  \makebox(0,0){\strut{}$\xi$\vspace{-7mm}}%
}}%
  \put(5155,440){\makebox(0,0){\strut{} 0.3}}%
  \put(4328,440){\makebox(0,0){\strut{} 0.1}}%
  \put(3421,440){\makebox(0,0){\strut{} 0.03}}%
  \put(2594,440){\makebox(0,0){\strut{} 0.01}}%
  \put(1687,440){\makebox(0,0){\strut{} 0.003}}%
  \put(860,440){\makebox(0,0){\strut{} 0.001}}%
  \put(740,2884){\makebox(0,0)[r]{\strut{}5.5}}%
  \put(740,2660){\makebox(0,0)[r]{\strut{}5.0}}%
  \put(740,2435){\makebox(0,0)[r]{\strut{}4.5}}%
  \put(740,2211){\makebox(0,0)[r]{\strut{}4.0}}%
  \put(740,1986){\makebox(0,0)[r]{\strut{}3.5}}%
  \put(740,1762){\makebox(0,0)[r]{\strut{}3.0}}%
  \put(740,1538){\makebox(0,0)[r]{\strut{}2.5}}%
  \put(740,1313){\makebox(0,0)[r]{\strut{}2.0}}%
  \put(740,1089){\makebox(0,0)[r]{\strut{}1.5}}%
  \put(740,864){\makebox(0,0)[r]{\strut{}1.0}}%
  \put(740,640){\makebox(0,0)[r]{\strut{}0.5}}%
\includegraphics{cbpdn_newexp02_lmdxi_64}%
\end{picture}%
\endgroup
 

%% file: cbpdn_newexp02_mnitlmbd_64.pstex_t
% GNUPLOT: LaTeX picture with Postscript
\begingroup%
\makeatletter%
\newcommand{\GNUPLOTspecial}{%
  \@sanitize\catcode`\%=14\relax\special}%
\setlength{\unitlength}{0.0500bp}%
\begin{picture}(5615,3124)(0,0)%
  \put(3386,1934){\makebox(0,0)[r]{\strut{}$\xi = 2.0$}}%
  \put(3386,2134){\makebox(0,0)[r]{\strut{}$\xi = 1.0$}}%
  \put(3386,2334){\makebox(0,0)[r]{\strut{}$\xi = 0.5$}}%
  \put(3386,2534){\makebox(0,0)[r]{\strut{}$\xi = f(\lambda)$}}%
  \put(3386,2734){\makebox(0,0)[r]{\strut{}Best $\xi$}}%
  \put(3007,140){\makebox(0,0){\strut{}$\lambda$}}%
  \put(160,1762){%
\rotatebox{-270}{%
  \makebox(0,0){\strut{}Mean iterations\vspace{-7mm}}%
}}%
  \put(5155,440){\makebox(0,0){\strut{} 0.3}}%
  \put(4328,440){\makebox(0,0){\strut{} 0.1}}%
  \put(3421,440){\makebox(0,0){\strut{} 0.03}}%
  \put(2594,440){\makebox(0,0){\strut{} 0.01}}%
  \put(1687,440){\makebox(0,0){\strut{} 0.003}}%
  \put(860,440){\makebox(0,0){\strut{} 0.001}}%
  \put(740,2884){\makebox(0,0)[r]{\strut{}500}}%
  \put(740,2635){\makebox(0,0)[r]{\strut{}450}}%
  \put(740,2385){\makebox(0,0)[r]{\strut{}400}}%
  \put(740,2136){\makebox(0,0)[r]{\strut{}350}}%
  \put(740,1887){\makebox(0,0)[r]{\strut{}300}}%
  \put(740,1637){\makebox(0,0)[r]{\strut{}250}}%
  \put(740,1388){\makebox(0,0)[r]{\strut{}200}}%
  \put(740,1139){\makebox(0,0)[r]{\strut{}150}}%
  \put(740,889){\makebox(0,0)[r]{\strut{}100}}%
  \put(740,640){\makebox(0,0)[r]{\strut{} 50}}%
\includegraphics{cbpdn_newexp02_mnitlmbd_64}%
\end{picture}%
\endgroup
 

%% file: admm.bbl
\begin{thebibliography}{10}
\providecommand{\url}[1]{#1}
\csname url@samestyle\endcsname
\providecommand{\newblock}{\relax}
\providecommand{\bibinfo}[2]{#2}
\providecommand{\BIBentrySTDinterwordspacing}{\spaceskip=0pt\relax}
\providecommand{\BIBentryALTinterwordstretchfactor}{4}
\providecommand{\BIBentryALTinterwordspacing}{\spaceskip=\fontdimen2\font plus
\BIBentryALTinterwordstretchfactor\fontdimen3\font minus
  \fontdimen4\font\relax}
\providecommand{\BIBforeignlanguage}[2]{{%
\expandafter\ifx\csname l@#1\endcsname\relax
\typeout{** WARNING: IEEEtran.bst: No hyphenation pattern has been}%
\typeout{** loaded for the language `#1'. Using the pattern for}%
\typeout{** the default language instead.}%
\else
\language=\csname l@#1\endcsname
\fi
#2}}
\providecommand{\BIBdecl}{\relax}
\BIBdecl

\bibitem{goldstein-2009-split}
T.~Goldstein and S.~J. Osher, ``The split {B}regman method for l1-regularized
  problems,'' \emph{SIAM Journal on Imaging Sciences}, vol.~2, no.~2, pp.
  323--343, 2009.  \doi{10.1137/080725891}

\bibitem{boyd-2010-distributed}
S.~Boyd, N.~Parikh, E.~Chu, B.~Peleato, and J.~Eckstein, ``Distributed
  optimization and statistical learning via the alternating direction method of
  multipliers,'' \emph{Foundations and Trends in Machine Learning}, vol.~3,
  no.~1, pp. 1--122, 2010.  \doi{10.1561/2200000016}

\bibitem{afonso-2011-augmented}
M.~V. Afonso, J.~M. Bioucas-Dias, and M.~A.~T. Figueiredo, ``An {A}ugmented
  {L}agrangian approach to the constrained optimization formulation of imaging
  inverse problems,'' \emph{IEEE Transactions on Image Processing}, vol.~20,
  no.~3, pp. 681--695, Mar. 2011.  \doi{10.1109/tip.2010.2076294}

\bibitem{ghadimi-2015-optimal}
E.~Ghadimi, A.~Teixeira, I.~Shames, and M.~Johansson, ``Optimal parameter
  selection for the alternating direction method of multipliers ({ADMM}):
  Quadratic problems,'' \emph{IEEE Transactions on Automatic Control}, vol.~60,
  no.~3, pp. 644--658, Mar. 2015.  \doi{10.1109/TAC.2014.2354892}

\bibitem{raghunathan-2014-alternating}
A.~U. Raghunathan and S.~{Di Cairano}, ``Alternating direction method of
  multipliers for strictly convex quadratic programs: optimal parameter
  selection,'' in \emph{American Control Conference (ACC)}, Jun. 2014, pp.
  4324--4329.  \doi{10.1109/ACC.2014.6859093}

\bibitem{raghunathan-2015-admm}
------, ``{ADMM} for convex quadratic programs: Linear convergence and
  infeasibility detection,'' arXiv, Tech. Rep. arXiv:1411.7288v2, 2015.

\bibitem{he-2000-alternating}
B.-S. He, H.~Yang, and S.-L. Wang, ``Alternating direction method with
  self-adaptive penalty parameters for monotone variational inequalities,''
  \emph{Journal of Optimization Theory and Applications}, vol. 106, pp.
  337--356, 2000.  \doi{10.1023/a:1004603514434}

\bibitem{hansson-2012-subspace}
\BIBentryALTinterwordspacing
A.~Hansson, Z.~Liu, and L.~Vandenberghe, ``Subspace system identification via
  weighted nuclear norm optimization,'' \emph{CoRR}, vol. abs/1207.0023, 2012.
  [Online]. Available: \url{http://arxiv.org/abs/1207.0023}
\BIBentrySTDinterwordspacing

\bibitem{liu-2013-nuclear}
Z.~Liu, A.~Hansson, and L.~Vandenberghe, ``Nuclear norm system identification
  with missing inputs and outputs,'' \emph{Systems \& Control Letters},
  vol.~62, no.~8, pp. 605 -- 612, 2013.  \doi{10.1016/j.sysconle.2013.04.005}

\bibitem{vu-2013-fantope}
V.~Q. Vu, J.~Cho, J.~Lei, and K.~Rohe, ``Fantope projection and selection: A
  near-optimal convex relaxation of sparse {PCA},'' in \emph{Advances in Neural
  Information Processing Systems 26}, C.~J.~C. Burges, L.~Bottou, M.~Welling,
  Z.~Ghahramani, and K.~Q. Weinberger, Eds., 2013, pp. 2670--2678.

\bibitem{iordache-2014-collaborative}
M.-D. Iordache, J.~M. Bioucas-Dias, and A.~Plaza, ``Collaborative sparse
  regression for hyperspectral unmixing,'' \emph{IEEE Transactions on
  Geoscience and Remote Sensing}, vol.~52, no.~1, pp. 341--354, Jan. 2014.
  \doi{10.1109/TGRS.2013.2240001}

\bibitem{weller-2014-phase}
D.~S. Weller, A.~Pnueli, O.~Radzyner, G.~Divon, Y.~C. Eldar, and J.~A. Fessler,
  ``Phase retrieval of sparse signals using optimization transfer and {ADMM},''
  \emph{Proc. IEEE Intl. Conf. on Image Processing}, pp. 1342--6, 2014.

\bibitem{wohlberg-2014-efficient}
B.~Wohlberg, ``Efficient convolutional sparse coding,'' in \emph{Proceedings of
  IEEE International Conference on Acoustics, Speech, and Signal Processing
  (ICASSP)}, Florence, Italy, May 2014, pp. 7173--7177.
  \doi{10.1109/ICASSP.2014.6854992}

\bibitem{urruty-2004-fundamentals}
J.-B. Hiriart-Urruty and C.~Lemar\'echal, \emph{Fundamentals of Convex
  Analysis}.\hskip 1em plus 0.5em minus 0.4em\relax Springer, 2004.

\bibitem{eckstein-2012-augmented}
\BIBentryALTinterwordspacing
J.~Eckstein, ``{A}ugmented {L}agrangian and alternating direction methods for
  convex optimization: A tutorial and some illustrative computational
  results,'' Rutgers Center for Operations Research, Rutgers University, Rutcor
  Research Report RRR 32-2012, December 2012. [Online]. Available:
  \url{http://rutcor.rutgers.edu/pub/rrr/reports2012/32_2012.pdf}
\BIBentrySTDinterwordspacing

\bibitem{wang-2001-decomposition}
S.-L. Wang and L.~Z. Liao, ``Decomposition method with a variable parameter for
  a class of monotone variational inequality problems,'' \emph{Journal of
  Optimization Theory and Applications}, vol. 109, pp. 415--429, 2001.
  \doi{10.1023/a:1017522623963}

\bibitem{mittelmann-2003-independent}
H.~D. Mittelmann, ``\BIBforeignlanguage{English}{An independent benchmarking of
  {SDP} and {SOCP} solvers},'' \emph{\BIBforeignlanguage{English}{Mathematical
  Programming}}, vol.~95, no.~2, pp. 407--430, 2003.
  \doi{10.1007/s10107-002-0355-5}

\bibitem{wachter-2006-implementation}
A.~W\"achter and L.~T. Biegler, ``\BIBforeignlanguage{English}{On the
  implementation of an interior-point filter line-search algorithm for
  large-scale nonlinear programming},''
  \emph{\BIBforeignlanguage{English}{Mathematical Programming}}, vol. 106,
  no.~1, pp. 25--57, 2006.  \doi{10.1007/s10107-004-0559-y}

\bibitem{ramdas-2015-fast}
A.~Ramdas and R.~J. Tibshirani, ``Fast and flexible admm algorithms for trend
  filtering,'' \emph{Journal of Computational and Graphical Statistics},
  vol.~25, no.~3, pp. 839--858, 2016.  \doi{10.1080/10618600.2015.1054033}

\bibitem{chen-1998-atomic}
S.~S. Chen, D.~L. Donoho, and M.~A. Saunders, ``Atomic decomposition by basis
  pursuit,'' \emph{SIAM J. Sci. Comput.}, vol.~20, no.~1, pp. 33--61, 1998.
  \doi{10.1137/S1064827596304010}

\bibitem{wohlberg-2016-efficient}
B.~Wohlberg, ``Efficient algorithms for convolutional sparse representations,''
  \emph{IEEE Transactions on Image Processing}, vol.~25, no.~1, pp. 301--315,
  Jan. 2016.  \doi{10.1109/TIP.2015.2495260}

\bibitem{wohlberg-2016-sporco}
------, ``{SP}arse {O}ptimization {R}esearch {CO}de ({SPORCO}),'' Software
  library available from \url{http://purl.org/brendt/software/sporco}, 2016.

\bibitem{parikh-2014-block}
N.~Parikh and S.~Boyd, ``Block splitting for distributed optimization,''
  \emph{Mathematical Programming Computation}, vol.~6, no.~1, pp. 77--102,
  2014.  \doi{10.1007/s12532-013-0061-8}

\bibitem{dyer-2013-greedy}
\BIBentryALTinterwordspacing
E.~L. Dyer, A.~C. Sankaranarayanan, and R.~G. Baraniuk, ``Greedy feature
  selection for subspace clustering,'' \emph{Journal of Machine Learning
  Research}, vol.~14, pp. 2487--2517, 2013. [Online]. Available:
  \url{http://jmlr.org/papers/v14/dyer13a.html}
\BIBentrySTDinterwordspacing

\bibitem{alliney-1992-digital}
S.~Alliney, ``Digital filters as absolute norm regularizers,'' \emph{IEEE
  Transactions on Signal Processing}, vol.~40, no.~6, pp. 1548--1562, Jun.
  1992.  \doi{10.1109/78.139258}

\bibitem{esser-2010-primal}
E.~Esser, ``Primal dual algorithms for convex models and applications to image
  restoration, registration and nonlocal inpainting,'' Ph.D. dissertation,
  University of California, Los Angeles, 2010.

\end{thebibliography}
